%% file: graph10.tex
\documentclass{amsart}
\usepackage{bm}
\usepackage{amsmath,amssymb}
\usepackage{graphicx}

\newtheorem{theorem}{Theorem}[section]
\newtheorem{lemma}[theorem]{Lemma}
\newtheorem{prop}[theorem]{Proposition}
\newtheorem{cor}[theorem]{Corollary}

\theoremstyle{definition}
\newtheorem{definition}[theorem]{Definition}
\newtheorem{example}[theorem]{Example}

\theoremstyle{remark}
\newtheorem{remark}[theorem]{Remark}

\begin{document}
\newcommand{\halfintegers}{\mathbb{Z}+\frac{1}{2}}
\newcommand{\aut}{{\rm aut}}
\input{graphexp2.tex}

\title{Pole Structure of Topological String Free Energy}
\author{Yukiko Konishi}

\thanks{Communicated by K.Saito. Received November 24, 2004.
Revised March 31, 2005.}
\thanks{
\textup{2000} {\it Mathematics Subject Classification.} 
{Primary 14N35; Secondary 05E05.}}

\thanks{Research Institute for Mathematical Sciences,
Kyoto University, Kyoto 606-8502, Japan.}
\thanks{
E-mail address: konishi@kurims.kyoto-u.ac.jp}

\begin{abstract}
We show that the free energy of the topological string
admits a certain pole structure by using the operator
formalism. Combined with the results of Peng that proved the
integrality,
this gives a combinatoric proof of the Gopakumar--Vafa conjecture.
\end{abstract}

\maketitle

\section{Introduction and Main Theorem} %
Recent developments in the string duality 
make it possible to express the 
partition function and the free energy 
of the topological string on
a toric Calabi--Yau threefold
in terms of the symmetric functions 
(see \cite{akmv}).
In 
mathematical terms,
the free energy is none other than
the generating function of
the Gromov--Witten invariants
%
%
\cite{lllz}.
In this paper we treat the case of 
the canonical
bundle of a toric surface
and 
its straightforward generalization.

The partition function is given as follows.
Let $r\geq 2$ be an integer  and
$\gamma=(\gamma_1,\ldots,\gamma_r)$ be
an $r$-tuple of integers  which will be fixed
from here on.
Let $\Vec{Q}=(Q_1,\ldots,Q_r)$ be an 
$r$-tuple of (formal) variables
and $q$ a variable.
%
\begin{definition}
\label{part}%
\begin{equation}\notag
\begin{split}
&{\mathcal Z}^{\gamma}(q;\Vec{Q})
=1+
\sum_{\begin{subarray}{l}\Vec{d}\in \mathbb{Z}^r_{\geq0},
\Vec{d}\neq \Vec{0}\end{subarray}}
{\mathcal Z}^{\gamma}_{\Vec{d}}(q) \, \Vec{Q}^{\Vec{d}},
\\
&{\mathcal Z}^{\gamma}_{\Vec{d}}(q)=
(-1)^{\gamma\cdot\Vec{d}}
\sum_{\begin{subarray}{c}(\lambda^1,\ldots,\lambda^r)\\
\lambda^i\in{\mathcal P_{d_i}}
\end{subarray}}
\prod_{i=1}^rq^{\frac{\gamma_i\kappa(\lambda^i)}{2}}
W_{\lambda^i,\lambda^{i+1}}(q),
\end{split}
\end{equation}
where 
$\Vec{Q}^{\Vec{d}}={Q_1}^{d_1}\cdots {Q_r}^{d_r}$
for $\Vec{d}=(d_1,\ldots,d_r)$,
${\mathcal P}_d$ is the set of partitions of $d$,
$\lambda^{r+1}=\lambda^1$ and
\begin{equation}\notag
\kappa(\lambda)=\sum_{i\geq 1}\lambda_i(\lambda_i-2i+1)
\qquad (\lambda\text{: a partition}).
\end{equation}
For a pair of partitions $(\mu,\nu)$,
$W_{\mu,\nu}(q)$ is defined as follows.
\begin{equation}\notag
W_{\mu,\nu}(q)=
(-1)^{|\mu|+|\nu|}q^{\frac{\kappa(\mu)+\kappa(\nu)}{2}}
\sum_{\eta\in {\mathcal P}}
s_{\mu/\eta}(q^{-\rho})s_{\nu/\eta}(q^{-\rho})
\end{equation}
where 
 ${\mathcal P}$ is the set of partitions and
$s_{\mu/\eta}(q^{-\rho})$ is the skew-Schur function 
associated to partitions $\mu,\eta$
with variables specialized to 
$q^{-\rho}=(q^{i-\frac{1}{2}})_{i\geq 1}$.
\end{definition}
We also define the free energy 
to be the logarithm of the partition function:
\begin{equation}\notag
{\mathcal F}^{\gamma}(q;\Vec{Q})=
\log {\mathcal Z}^{\gamma}(q;\Vec{Q}).
\end{equation}
The logarithm  should be considered as
a formal power series in the variables $Q_1,\ldots,Q_r$.
The coefficient of $\Vec{Q}^{\Vec{d}}$ is denoted by
${\mathcal F}_{\Vec{d}}^{\gamma}(q)$.

The free energy is related to
the Gromov--Witten invariants 
in the following way.
The target manifold is the total space $X$ of 
the canonical bundle on a smooth, complete toric surface $S$.
Recall that a  toric surface $S$ is 
given by a two-dimensional complete fan;
its one-dimensional cones (1-cones)
correspond to toric invariant
rational curves.
Let $r_S$ be the number of 1-cones,
$\Vec{C}=(C_1,\ldots,C_{r_S})$ the set of 
the rational curves 
and
${\gamma}_S$ 
the set of the self-intersection numbers:
\begin{equation}\notag
\gamma_{S}=(C_1^2,\ldots,C_r^2).
\end{equation}
For example,
if $S$ is 
${\bf P}^2$,
${\bf F}_0$, 
${\bf F}_1,{\mathcal B}_2$ or ${\mathcal B}_3$,
$\gamma_S$ is 
$(1,1,1)$,
$(0,0,0,0)$,
$(1,0,-1,0)$,
$(0,0,-1,-1,-1)$ or
$(-1,-1,-1,-1,-1,-1)$.
Then the generating function of 
the Gromov--Witten invariants $N_{\beta}^g(X)$ of $X=K_S$ 
with fixed degree $\beta$,
\begin{equation}\notag
{\mathcal F}_{\beta}(X)=\sum_{g\geq 0}N_{\beta}^g(X)g_s^{2g-2},
\end{equation}
is exactly equal to the sum \cite{z1}\cite{llz2}:
\begin{equation}\notag
{\mathcal F}_{\beta}(X)=\sum_{\Vec{d};[\Vec{d}\cdot\Vec{C}]=\beta}
{\mathcal F}_{\Vec{d}}^{\gamma_S}(q)\big|_{q=e^{\sqrt{-1}g_s}}.
\end{equation}
Actually, in the localization calculation,
each 
$ {\mathcal F}_{\Vec{d}}^{\gamma_S}(q)\big|_{q=e^{\sqrt{-1}g_s}}$
is the contribution from the fixed point loci
in the moduli of stable maps
of which the image curves are $\Vec{d}.\Vec{C}$
(see \cite{z1}).

Note that some
pairs $(r,\gamma)$ 
do not correspond to toric surfaces.
One of the simplest cases is $r=2$, since
for any two dimensional fan to be complete, it
must has 
at least three 1-cones.
In this article, 
we also deal with such non-geometric cases.

One problem concerning the Gromov--Witten invariants 
is the Gopakumar--Vafa conjecture
\cite{gv}. (see also \cite{bp}\cite{hst}).
Let us 
define
the numbers
$\{n_{\beta}^g(X)\}_{g,\beta}$ by rewriting
$\{{\mathcal F}_{\beta}(X)\}_{\beta\in H_2(X;\mathbb{Z})}$
in the form below.
\begin{equation}\label{kakinaosi}
{\mathcal F}_{\beta}(X)=\sum_{g\geq 0}\sum_{k;k|\beta}
\frac{n_{\beta/k}^g(X)}{k}
\Big(2\sin\frac{kg_s}{2}\Big)^{2g-2}.
\end{equation}
Then the conjecture states the followings.
\begin{enumerate}
\item\label{integrality}
$n_{\beta}^g(X)\in\mathbb{Z}$ and
$n_{\beta}^g(X)=0$ for every fixed $\beta$ and $g\gg1$.
\item
\label{bps}
Moreover,
$n_{\beta}^g(X)$ is equal to the 
number of certain BPS states in M-theory
(see \cite{hst} for a mathematical formulation).
\end{enumerate}
By 
the M\"obius inversion formula \cite{bp},
(\ref{kakinaosi}) is equivalent to
\begin{equation}\label{sononi}
\sum_{k;k|\beta}\frac{\mu(k)}{k}
{\mathcal F}_{\beta/k}(X)\big|_{g_s\to kg_s}
=
\sum_{g\geq 1}n_{\beta}^g(X)
\Big(2\sin\frac{g_s}{2}\Big)^{2g-2}
\end{equation}
where $\mu(k)$
is the M\"obius function.
Therefore the first part
is equivalent to the LHS being a
polynomial in $t=-\big(2\sin\frac{g_s}{2}\big)^2$ with integer
coefficients divided by $t$.

In this article, 
we deal with the first part of the Gopakumar--Vafa
conjecture.
As it turns out, it holds
not only for a class $\beta\in H_2(X;\mathbb{Z})$
but also for each
torus equivariant class $\Vec{d}\cdot\Vec{C}
$.
Moreover, it also holds in non-geometric cases.

\begin{definition}%
We define
\begin{equation}\notag
G_{\Vec{ d}}^{\gamma}(q)=
\sum_{k';k'|k}\frac{k'}{k}\mu\Big(\frac{k}{k'}\Big)
{\mathcal F}_{k'{\Vec{ d}}/k}^{\gamma}(q^{k/k'})\qquad
(k={\rm gcd}(\Vec{d})).
\end{equation}
\end{definition} %
The main result of the paper is the following theorem.
We set $t=(q^{\frac{1}{2}}-q^{-\frac{1}{2}})^2$.
\begin{theorem} \label{mainthm}
$t \cdot G_{\Vec{d}}^{\gamma}(q)\in\mathbb{Q}[t]$.
\end{theorem}%
Peng proved
that $G_{\Vec{d}}^{\gamma}(q)$
is a rational function in $t$
such that its numerator and denominator
are polynomials with integer coefficients and 
the denominator is monic \cite{peng}.
So we have
\begin{cor}\label{maincor}
$t\cdot G_{\Vec{ d}}^{\gamma}(q) \in\mathbb{Z}[t]$.
\end{cor}%
In the geometric case corresponding to $X=K_S$,
the corollary implies the
integrality of $n_{\beta}^g(X)$ and its vanishing at
higher genera.
This is because
the LHS of (\ref{sononi}) is equal to
\begin{equation}\notag
\sum_{\Vec{d};[\Vec{d}\cdot\Vec{C}]=\beta}
G_{\Vec{d}}^{\gamma}(q)\big|_{q=e^{\sqrt{-1}g_s}}.
\end{equation}

The organization of the paper is as follows.
In section \ref{pfof}, we 
introduce the infinite-wedge space (Fock space)
and the fermion operator algebra
and write the partition function
in terms of matrix elements of a certain operator.
In section \ref{me}, we
express the matrix elements as the sum of certain 
quantities - amplitudes -  over a set of graphs.
Then we rewrite
the partition function
in terms of graph amplitudes.
In section \ref{feen}, 
we take the logarithm of the partition function and 
obtain the free energy.
The key idea is to
use the exponential formula,
which is
the relation between log/exp and connected/disconnected
graph sums.
We give an outline of the proof of the main theorem
in section \ref{pt}.
Then, in section \ref{ps}, we study
the pole structure of the amplitudes
and finish the proof.
The rigorous formulation of the exponential formula
and the proof of the free energy appear in
appendices \ref{exf} and \ref{pfe}.
Appendix \ref{proofG4}
contains the proof of a lemma. 

\begin{center}{\bf Acknowledgement}\end{center} %

The author thanks
H. Iritani, H. Kajiura, T. Kawai, K. Saito,
A. Takahashi and K. Ueda
for valuable discussions
and comments.
She also thanks M. Yoshinaga for recommending the book of Stanley
\cite{stan}.

The author thanks  S. Hosono
for
carefully reading the manuscript and 
giving useful suggestions for improvements.
%

\section{Partition Function in Operator Formalism}\label{pfof}

The goal of this section is to express the partition function
in terms of the matrix elements of certain operator in
the fermion operator algebra.
We first introduce the infinite wedge space \cite{op1} 
(also called the Fock space in \cite{knty})
and an action of the fermion operator algebra in subsection \ref{OP}.
Then we see that 
the skew-Schur function and other quantities of partitions
($|\lambda|$ and $\kappa(\lambda)$) are the matrix elements 
of operators.
In subsection \ref{PFOP},
we rewrite $W_{\mu,\nu}(q)$ and the partition function.
%
\subsection{Operator Formalism}\label{OP}
In this subsection 
we first briefly explain notations (mainly) on partitions.
%
Secondly we  
introduce the infinite wedge space,
the fermion operator algebra and define
some operators.
Then we restrict ourselves to
a subspace of the infinite wedge space.
We see that
the canonical basis is naturally associated to the set of
partitions and 
that 
the skew-Schur function, $|\lambda|$ and $\kappa(\lambda)$ are
the matrix elements of certain operators with respect to the basis.
Finally we introduce a new basis which will play an
important role in the later calculations.

\subsubsection{Partitions}\label{NP}
%

A {\em partition} is  a non-increasing sequence 
$\lambda=(\lambda_1,\lambda_2,\ldots)$ of 
non-negative integers containing only finitely many nonzero
terms.
The nonzero $\lambda_i$'s  are called the parts.
The number of parts is the {\em length} of $\lambda$, 
denoted by $l(\lambda)$. 
The sum of the parts is the {\em weight} of 
$\lambda$, denoted by $|\lambda|$: $|\lambda|=\sum_i\lambda_i$.
If $|\lambda|=d$, $\lambda$ is a partition of $d$. 
The set of all partitions of $d$ is denoted by ${\mathcal P}_d$
and the set of all partitions by ${\mathcal P}$. 
Let $m_k(\lambda)=\#\{\lambda_i:\lambda_i=k\}$ 
be the {\em multiplicity}
of $k$ where
$\#$ denotes the number of elements of a finite set. 
Let $\aut(\lambda)$ be the symmetric group 
acting as the permutations of  the equal parts of $\lambda$:
$\aut(\lambda)\cong \prod_{k\geq 1}\mathfrak{S}_{m_k(\lambda)}$.
Then $\#\aut(\lambda)=\prod_{k\geq 1}m_k(\lambda)!$.
We define
\begin{equation}\notag
z_{\lambda}=\prod_{i=1}^{l(\lambda)}\lambda_i \cdot\#\aut(\lambda),
\end{equation}
which is the number of the centralizers of the conjugacy class 
associated to $\lambda$.

A partition $\lambda=(\lambda_1,\lambda_2,\ldots)$ is identified as 
the Young diagram with $\lambda_i$ boxes in the $i$-th row 
$(1\leq i\leq l(\lambda))$.
The Young diagram with $\lambda_i$ boxes in the $i$-th column is 
its {\em transposed} Young diagram. The corresponding partition
is called the {\em conjugate partition} and denoted by $\lambda^t$.
Note that
$\lambda^t_i=\sum_{k\geq i}m_k(\lambda)$.

We define
\begin{equation}\notag
\kappa(\lambda)=\sum_{i=1}^{l(\lambda)}\lambda_i(\lambda_i-2i+1).
\end{equation}
This is equal to twice the sum of contents $\sum_{x\in \lambda}c(x)$ 
where $c(x)=j-i$ for the box $x$ at the $(i,j)$-th place 
in the Young diagram $\lambda$.
Thus, $\kappa(\lambda)$ is always even and 
satisfies $\kappa(\lambda^t)=-\kappa(\lambda)$.

$\mu\cup\nu$ denotes the partition whose parts 
are $\mu_1,\ldots,\mu_{l(\mu)}, 
\nu_1,\ldots,\nu_{l(\mu)}$ and 
$k\mu$ the partition $(k\mu_1,k\mu_2,\ldots)$ 
for $k\in {\mathbb{N}}$.

We define
\begin{equation}\notag
[k]=q^{\frac{k}{2}}-q^{-\frac{k}{2}} \qquad (k\in\mathbb{Q}),
\end{equation}
which is called the {\em $q$-number}.
For a partition $\lambda$, we use the shorthand notation
\begin{equation}\notag
[\lambda]=\prod_{i=1}^{l(\lambda)}[\lambda_i].
\end{equation}
\subsubsection{Infinite Wedge Space (Fock Space)}
%
A subset $S$ of $\halfintegers$ is called a {\em Maya diagram}
if both $S_+:=S\cap \{k\in\halfintegers|k>0\}$ and
$S_-:=S^c\cap\{k\in\halfintegers|k<0\}$ are finite sets.
The {\em charge} $\chi(S)$  is defined by
\begin{equation}\notag
\chi(S)=\# S_+-\#S_-.
\end{equation}
The set of all Maya diagrams of charge $p$ is
denoted by ${\mathcal M}_p$.
We write a Maya diagram $S$ in the decreasing sequence
$S=(s_1,s_2,s_3,\ldots)$. Note that if $\chi(S)=p$,
$s_i\geq p-i+\frac{1}{2}$
for all $i\geq 1$
and $s_i= p-i+\frac{1}{2}$
for $i\gg 1$.
Therefore
\begin{equation}\notag
\lambda^{(p)}(S)=\Big(s_1-p+\frac{1}{2},s_2-p+\frac{3}{2},\ldots\Big)
\end{equation}
is a partition. 
$\lambda^{(p)}:{\mathcal M}_p\stackrel{\sim}{\to}{\mathcal P}$
is a canonical bijection for each $p$
where the inverse 
$(\lambda^{(p)})^{-1}(\mu)=S$ $(\mu\in\mathcal P)$ 
is given by 
\begin{equation}\notag
S_+=\Big\{\mu_i-i+\frac{1}{2}+p\Big|1\leq i\leq k\Big\},
\qquad
S_-=\Big\{\mu_i^t-i+\frac{1}{2}+p\Big|1\leq i\leq k\Big\}.
\end{equation}
Here $k=\#\{\mu_i|\mu_i=i\}$ is the number of diagonal boxes
in the Young diagram of $\mu$.
We define 
\begin{equation}\notag
d(S)=|\lambda^{(p)}(S)|.
\end{equation}

Let $V$ be  an infinite dimensional linear space over ${\mathbb{C}}$
equipped with 
a basis $\{e_k\}_{k\in \halfintegers}$
satisfying the following condition:
every element $v\in V$ is expressed as
$v=\sum_{k>m}v_k e_k$ with some $m\in {\mathbb Z}$.
Let $\bar{V}={\rm Hom}_{\mathbb{C}}(V,{\mathbb{C}})$ be 
the topological dual space and $\{\bar{e}_k\}_{k\in\halfintegers}$ 
the dual basis:
$\bar{e}_l(e_k)=\delta_{k,l}$.

For each Maya diagram $S=(s_1,s_2,s_3,\ldots)$, 
 $|v_S\rangle$ denotes the symbol
\begin{equation}\notag
|v_S\rangle=e_{s_1}\wedge e_{s_2}\wedge e_{s_3}\wedge\cdots.
\end{equation}
The infinite wedge space of charge $p$, $\Lambda^{\frac{\infty}{2}}_pV$
is the vector space over ${\mathbb{C}}$
spanned by $\{|v_S\rangle\}_{S\in {\mathcal M}_p}$, 
and the infinite wedge space $\Lambda^{\frac{\infty}{2}}V$
is the direct sum of the charge $p$ spaces $(p\in\mathbb{Z})$:
\begin{equation}\notag
\Lambda^{\frac{\infty}{2}}_pV=\prod_{S\in{\mathcal M}_p}
\mathbb{C}|v_S\rangle,
\qquad 
\Lambda^{\frac{\infty}{2}}V=
\bigoplus_{p\in \mathbb{Z}}\Lambda^{\frac{\infty}{2}}_pV.
\end{equation}
%
We define
the charge operator $J_0$ and
the mass operator $M$ 
on $\Lambda^{\frac{\infty}{2}}V$ by
\begin{equation}\notag
J_0v_S=\chi(S) v_S,\qquad 
M v_S=d(S)v_S.
\end{equation}
%
Since $J_0$ and $M$ commute,  
$\Lambda^{\frac{\infty}{2}}_pV$  decomposes into the eigenspace 
$\Lambda^{\frac{\infty}{2}}_pV(d)$
of $M$
with eigenvalue $d$:
\begin{equation}\notag
\Lambda^{\frac{\infty}{2}}_pV(d)=
\bigoplus_{\chi(S)=p,d(S)=d}
\mathbb{C}|v_S\rangle.
\end{equation}
The dimension of the eigenspace is equal to 
the number $p(d)$ of
partitions of $d$.

For a Maya diagram $S=(s_1,s_2,\ldots)$, 
we define the symbol $\langle v_S|$ by
\begin{equation}
\notag
\langle v_S|=\ldots\wedge \bar{e}_{s_2}\wedge \bar{e}_{s_1}.
\end{equation}
The dual infinite wedge space is defined by
\begin{equation}\notag
\Lambda^{\frac{\infty}{2}}\bar{V}=
\bigoplus_{p\in\mathbb{Z}}
\Lambda^{\frac{\infty}{2}}_p\bar{V}
,\qquad
\Lambda^{\frac{\infty}{2}}_p\bar{V}
=
\bigoplus_{S\in {\mathcal M}_p}
\mathbb{C}\langle v_S|.
\end{equation}
The dual pairing is denoted by $\langle\phantom{y} |\phantom{y}\rangle$:
\begin{equation}\notag
\langle v_{S'}|v_S\rangle=\delta_{S,S'}.
\end{equation}

We set
\begin{equation}\notag
|p\rangle=e_{p-\frac{1}{2}}\wedge e_{p-\frac{3}{2}}\wedge\ldots,\qquad
\langle p|=\ldots \wedge \bar{e}_{p-\frac{3}{2}}
\wedge \bar{e}_{p-\frac{1}{2}}.
\end{equation}
$|p\rangle$ is called the {\em vacuum state} of the charge $p$.
Note that 
every $|p\rangle$ $(p\in\mathbb{Z})$ 
corresponds to the empty partition.
It is the basis of the subspace with charge $p$
and degree zero:
$\Lambda_p^{\frac{\infty}{2}}V(0)=\mathbb{C}|p\rangle$.

\subsubsection{Fermion Operator Algebra}%
Now we introduce the fermion operators algebra.
It is the associative algebra with $1$ generated by
$\psi_k,\psi^*_k$ $(k\in\halfintegers)$, with the relations:
\begin{equation}\notag
\{\psi_k,\psi_l^*\}=\delta_{k,l},\qquad
\{\psi_k,\psi_l\}=\{\psi_k^*,\psi_l^*\}=0
\end{equation}
for all $k,l\in\halfintegers$. 
Here $\{A,B\}=AB+BA$ is the anti-commutator.
We define the actions of $\psi_k,\psi_l^*$ on 
$\Lambda^{\frac{\infty}{2}}{V}$ and
$\Lambda^{\frac{\infty}{2}}\bar{V}$ as follows:
\begin{equation}
\notag
\begin{split}
\psi_k=e_k\wedge,\qquad \psi_k^*=\frac{\partial}{\partial e_k}
\qquad(\text{left action on $\Lambda^{\frac{\infty}{2}}{V}$}),\\
\psi_k=\frac{\partial}{\partial \bar{e}_{k}},\qquad 
\psi_k^*=\wedge \bar{e}_{k}
\qquad(\text{right action on $\Lambda^{\frac{\infty}{2}}\bar{V}$}).
\end{split}
\end{equation}
These are compatible with the dual pairing.
So any operator $A$ of the fermion algebra
satisfies
$\langle v_{S'}|(A|v_S\rangle)=
(\langle v_{S'}|A)|v_S\rangle$.
We call it the matrix element of $A$ with respect to 
$\langle v_{S'}|$ and $|v_{S}\rangle$ and write it as
\begin{equation}\notag
\langle v_{S'}|A|v_S\rangle.
\end{equation}
The operators $\psi_k$ and $\psi_k^*$ satisfy:
\begin{equation}\notag
\begin{split}
\psi_k|p\rangle =0 (k<p),\qquad
\psi_k^*|p\rangle =0(k>p),
\\
\langle p|\psi_k =0 (k>p),\qquad
\langle p|\psi^*_k=0 (k<p).
\end{split}
\end{equation}

Let us  define some operators. 
We first define, for $i,j\in\mathbb{Z}+\frac{1}{2},$
\begin{equation}\notag
E_{i,j}=:\psi_i\psi_j^*: ,
\qquad
:\psi_i\psi_j^*: =\begin{cases}
\psi_i\psi_j^*&(j>0)\\
-\psi_j^*\psi_i&(j<0).
\end{cases}
\end{equation}
The commutation relation is 
\begin{equation}\notag
[E_{i,j},E_{k,l}]=\delta_{j,k}E_{i,l}-\delta_{i,l}E_{k,j}+
\delta_{i,l}\delta_{j,k}(\theta(l<0)-\theta(j<0))
\end{equation}
where $\theta(l<0)=0$ if $l>0$ and $1$ if $l<0$.
$\theta(j<0)$ is defined similarly.
Next
we define
\begin{equation}\notag
C=\sum_{k\in{\mathbb{Z}}+\frac{1}{2}}E_{k,k},\qquad
H=\sum_{k\in{\mathbb{Z}}+\frac{1}{2}}kE_{k,k}
,\qquad
{\mathcal F}_2=\sum_{k\in{\mathbb{Z}}+\frac{1}{2}}\frac{k^2}{2}E_{k,k}.
\end{equation}
These act on the state $|v_S\rangle$ of charge $\chi(S)=p$  
as follows%
:
\begin{equation}\notag
\begin{split}
&C |v_S \rangle =p |v_S \rangle,\qquad
H |v_S \rangle =\Big(d(S)+\frac{p}{2}\Big)
|v_S \rangle,\\
&{\mathcal F}_2 |v_S \rangle =\Big(\frac{\kappa(\lambda^{(p)}(S))}{2}+
p\,d(S)+
\frac{p(4p^2-1)}{24} \Big)|v_S \rangle.
\end{split}
\end{equation}
Since $C$ is equal to the charge operator $J_0$,
it is also called the charge operator. 
We call $H$ the {\em energy} operator.

We define
\begin{equation}\notag
\alpha_m=\sum_{k\in\halfintegers}E_{k-m,k}
\qquad (m\in \mathbb{Z}\setminus\{0\}).
\end{equation}
Since these operators satisfy the commutation relations 
$[\alpha_m,\alpha_n]=m\delta_{m+n,0}$, they are called {\em bosons}.
Note that
$[C,\alpha_m]=0$, $[H,\alpha_m]=-m\alpha_m$.
%

The operator
\begin{equation}\notag
\Gamma_{\pm}(p)=\exp\Big[\sum_{{n\geq 1}}\frac{p_n\alpha_{\pm n}}{n}
\Big].
\end{equation}
is called the {\em vertex operators}
where $p=(p_1,p_2,\ldots)$ is a (possibly infinite) sequence.
In the later calculation, the sequence $p$ is taken to be
the power sum functions of certain variables.

Next we define the operator (see \cite{op1})
which will play an important role later.
\begin{equation}\notag
{\mathcal E}_{c}(n)=
\sum_{k\in{\mathbb{Z}}+\frac{1}{2}}q^{n(k-\frac{c}{2})}E_{k-c,k}
+\frac{\delta_{c,0}}{[n]}
\qquad \big((c,n)\in \mathbb{Z}^2\setminus\{(0,0)\}\big).
\end{equation}
This is, in a sense, a deformation of the boson $\alpha_m$ by 
${\mathcal F}_2$ since
\begin{equation}\label{f2alpha}
{\mathcal E}_m(0)=\alpha_m \quad (m\neq 0) 
\qquad
\text{ and }
\qquad
q^{{\mathcal F}_2}\alpha_{-n} q^{-{\mathcal F}_2}=
{\mathcal E}_{-n}(n) \quad (n\in \mathbb{N}).
\end{equation}
The commutation relation is as follows  
\begin{equation}\notag\label{commrel}
[{\mathcal E}_{a}(m),{\mathcal E}_{b}(n)]=
\begin{cases}
\Bigg[\,
\begin{vmatrix}a&m\\b&n\\\end{vmatrix}\,\Bigg]
{\mathcal E}_{a+b}(m+n)
&\text{ if }
(a+b,m+n)\neq (0,0)\\
a
&\text{ if }(a+b,m+n)=(0,0)
\end{cases}
\end{equation}
where $|\phantom{a}|$ in the RHS means the determinant
and $[\phantom{a}]$ the $q$-number.
Note that the operators $H,{\mathcal F}_2, \alpha_m,
\Gamma_{\pm}(p)$ and ${\mathcal E}_c(n)$ 
preserve the charge of a state 
because
they commute with the charge operator $C$.

\subsubsection{Charge Zero Subspace}
From here on  we  work only on the charge zero 
space $\Lambda_0^{\frac{\infty}{2}}V$.

The 
$\langle 0|A|0 \rangle$ is called
the {\em vacuum expectation value} (VEV) 
of the operator $A$ and denoted by
$
\langle A \rangle
$.

We use the set of partitions ${\mathcal P}$
to label the states instead of the set ${\mathcal M_0}$
of Maya diagrams: for a partition $\lambda$,
\begin{equation}\notag
|v_{\lambda}\rangle
=e_{\lambda_1-\frac{1}{2}}\wedge
e_{\lambda_2-\frac{3}{2}}\wedge\ldots,\qquad
\langle v_{\lambda}|=\ldots\wedge \bar{e}_{\lambda_2-\frac{3}{2}}
\wedge \bar{e}_{\lambda_1-\frac{1}{2}}.
\end{equation}
$|v_{\emptyset}\rangle=|0\rangle$, 
is the vacuum state and 
$\langle v_{\emptyset}|=\langle 0|$.

With this notation, the action of the energy operator $H$ 
and ${\mathcal F}_2$
are written as follows:
\begin{equation}\label{hk}
H |v_{\lambda} \rangle =|\lambda||v_{\lambda} \rangle,
\qquad
{\mathcal F}_2 |v_{\lambda} \rangle 
=\frac{\kappa(\lambda)}{2}|v_{\lambda} \rangle.
\end{equation}

One of the important points is 
that the Schur function and the 
skew Schur function are 
the matrix elements of the vertex operators:
let $x=(x_1,x_2,\ldots)$ be variables,
$p_i(x)=\sum_j (x_j)^i$ be the $i$-th power sum function
and $p(x)=(p_1(x),p_2(x),\ldots)$; 
then
\begin{equation}\label{opsym}
\begin{split}
&\langle v_{\lambda}|\Gamma_-(p(x))|0\rangle=
\langle 0|\Gamma_+(p(x))|v_{\lambda}\rangle=
s_{\lambda}(x),
\\
&\langle v_{\lambda}|\Gamma_-(p(x))|v_{\eta}\rangle=
\langle v_{\eta}|\Gamma_+(p(x))|v_{\lambda}\rangle=
s_{\lambda/\eta}(x).
\end{split}
\end{equation}


\subsubsection{Bosonic Basis}
For a partition $\mu\in{\mathcal P}$, we write
\begin{equation}\notag
\alpha_{\pm \mu}=\prod_{i=1}^{l(\mu)}\alpha_{\pm \mu_i}
.
\end{equation}
This is well-defined since
it is a product of commuting operators.
Now we
introduce the different kind of states:
\begin{equation}\notag
|\mu \rangle=\alpha_{-\mu}
|0\rangle,\qquad
\langle \mu|=\langle 0|\alpha_{\mu}
\qquad(\mu\in {\mathcal P}).
\end{equation}
The following relations are easily calculated
from the commutation relations:
\begin{equation}\notag
C|\mu\rangle=0,\qquad
H|\mu\rangle=|\mu||\mu\rangle
,\qquad
\langle\mu|\nu\rangle=z_{\mu}\delta_{\mu,\nu}.
\end{equation}
The set $\{|\mu\rangle |\mu\in{\mathcal P}_d\}$
is a basis of the degree $d$
subspace $\Lambda_0^{\frac{\infty}{2}}V(d)$
since it
consists of $p(d)=\dim \Lambda_0^{\frac{\infty}{2}}V(d)$ 
linearly independent states
of energy $d$.
\begin{equation}\notag
\Lambda_0^{\frac{\infty}{2}}V(d)=
\bigoplus_{\mu\in{\mathcal P}_d}\mathbb{C}|\mu\rangle.
\end{equation}
We call the new basis $\{|\mu\rangle| \mu\in {\mathcal P}_d\}$ 
the {\rm bosonic basis}
and the old basis $\{|v_\lambda\rangle| \lambda\in {\mathcal P}_d\}$ 
the fermionic basis.

The relation between the bosonic and fermionic basis
is written in terms of the character of the symmetric group.
Every irreducible, finite dimensional representation of
the symmetric group 
$\mathfrak{S}_{d}$
corresponds, one-to-one,
to a partition $\lambda\in{\mathcal P}_d$.
On the other hand, every conjugacy class 
also corresponds, one-to-one,
to a partition $\mu\in{\mathcal P}_d$.
The character of
the representation $\lambda$
on the 
conjugacy class ${\mu}$
is denoted by
$\chi_{\lambda}({\mu})$.
%
\begin{lemma}
\begin{equation}\notag
|v_{\lambda}\rangle
=\sum_{\mu\in{\mathcal P}_{|\lambda|}}
\frac{\chi_{\lambda}(\mu)}{z_{\mu}}
|\mu\rangle
\qquad
|\mu\rangle
=\sum_{\lambda\in{\mathcal P}_{|\mu|}}
\chi_{\lambda}(\mu)
|v_{\lambda}\rangle.
\end{equation}
\end{lemma}
\begin{proof}%
The vertex operator is formerly expanded as follows:
\begin{equation}\label{vtxop}
\Gamma_{\pm}(p(x))=\sum_{\mu\in{\mathcal P}}
\frac{p_{\mu}}{z_{\mu}}\alpha_{\pm\mu}.
\end{equation}
where
$p_{\mu}(x)=\prod_{i=1}^{l(\mu)}p_{\mu_i}(x)$.
Substituting into the first equation of (\ref{opsym}), we 
obtain
\begin{equation}\notag
s_{\lambda}(x)=\sum_{\mu\in{\mathcal P}_{|\lambda|}}
\langle v_{\lambda}|\mu\rangle \frac{p_{\mu}(x)}{z_{\mu}}
\end{equation}
Comparing with the Frobenius character formula
\begin{equation}\label{fro}
s_{\lambda}(x)= \sum_{\mu\in{\mathcal P}_{|\lambda|}}
\chi_{\lambda}(\mu) \frac{p_{\mu}(x)}{z_{\mu}},
\end{equation}
we obtain 
$\langle v_{\lambda}|\mu\rangle 
=\chi_{\lambda}(\mu).
$ 
Then the lemma follows immediately.
\end{proof}%
%
\subsection{Partition Function in Operator Formalism}\label{PFOP}
In this subsection, we first 
express $W_{\mu,\nu}(q)$ in terms of the fermion operator algebra.
Then we rewrite the partition function 
in terms of the matrix elements of the operator
$q^{{\mathcal F}_2}$ with respect to the bosonic basis.

\begin{lemma}\label{wmunu}
\begin{equation}\notag
W_{\mu,\nu}(q)=(-1)^{|\mu|+|\nu|}
\sum_{\eta';|\eta'|\leq |\mu|,|\nu|}
\sum_{\begin{subarray}{c}
\mu';|\mu'|=|\mu|-|\eta'|,\\
\nu';|\nu'|=|\nu|-|\eta'|
\end{subarray}}
\frac{(-1)^{l(\mu')+l(\nu')}}{z_{\mu'}z_{\nu'}z_{\eta'}[\mu'][\nu']}
\langle v_{\mu}|q^{{\mathcal F}_2}|\mu'\cup\eta'\rangle
\langle v_{\nu}|q^{{\mathcal F}_2}|\nu'\cup\eta'\rangle.
\end{equation}
\end{lemma}%
\begin{proof}%
We rewrite the skew-Schur function
in the variables $x=(x_1,x_2,\ldots)$.
By (\ref{opsym}) and (\ref{vtxop}),
\begin{equation}\notag
s_{\mu/\eta}(x)=
\sum_{\begin{subarray}{c}\mu';|\mu'|=|\mu|-|\eta|,\\
                         \eta';|\eta'|=|\eta|
      \end{subarray}}
\frac{p_{\mu'}(x)}{z_{\mu'}z_{\eta'}}
\langle v_{\mu}|\mu'\cup\eta'\rangle
\langle \eta'|v_{\eta}\rangle.
\end{equation}
When $\eta=\emptyset$, 
this is nothing but the Frobenius character formula (\ref{fro}).

The power sum function $p_i(x)$
associated to the specialized variables $q^{-\rho}$ is
\begin{equation}\notag p_i(q^{-\rho})=-\frac{1}{[\,i\,]}.\end{equation}
Therefore
\begin{equation}\notag
\begin{split}
&\sum_{\eta} s_{\mu/\eta}(q^{-\rho})s_{\nu/\eta}(q^{-\rho})\\
&=
\sum_{d=0}^{\min\{|\mu|,|\nu|\}}
\sum_{\eta\in{\mathcal P}_d}
s_{\mu/\eta}(q^{-\rho})s_{\nu/\eta}(q^{-\rho})
\\
&=\sum_{d=0}^{\min\{|\mu|,|\nu|\}}
\sum_{\begin{subarray}{c}
\mu';|\mu'|=|\mu|-d,\\
\nu';|\nu'|=|\nu|-d,\\
\lambda,\lambda'\in{\mathcal P}_d
\end{subarray}}
\frac{(-1)^{l(\mu')+l(\nu')}}
{z_{\mu'}z_{\nu'}z_{\lambda}z_{\lambda'}[\mu'][\nu']}
\langle v_{\mu}|\mu'\cup\lambda\rangle
\langle v_{\nu}|\nu'\cup\lambda'\rangle
\underbrace{
\sum_{\eta\in{\mathcal P}_d}
\langle\lambda|v_{\eta} \rangle
\langle v_{\eta}|\lambda' \rangle
}_{=\langle\lambda|\lambda'\rangle=z_{\lambda}\delta_{\lambda,\lambda'}}
\\
&=\sum_{\mu',\nu',\lambda}
\frac{(-1)^{l(\mu')+l(\nu')}}
{z_{\mu'}z_{\nu'}z_{\lambda}[\mu'][\nu']}
\langle v_{\mu}|\mu'\cup\lambda\rangle
\langle v_{\nu}|\nu'\cup\lambda\rangle.
\end{split}
\end{equation}

By (\ref{hk}), 
the factor $q^{\frac{\kappa(\mu)+\kappa(\nu)}{2}}$ is written
as follows:
\begin{equation}\notag
q^{\frac{\kappa(\mu)+\kappa(\nu)}{2}}=
\langle v_{\mu}|q^{{\mathcal F}_2}|v_{\mu}\rangle
\langle v_{\nu}|q^{{\mathcal F}_2}|v_{\nu}\rangle.
\end{equation}
Combining the above expressions, 
we obtain the lemma.
\end{proof}

Before rewriting the partition function, we explain some notations.
We use the symbols $\Vec{\mu},\Vec{\nu},\Vec{\lambda}$
to denote  $r$-tuples of partitions.
The $i$-th partition of $\Vec{\mu}$ 
is denoted by $\mu^i$
and $\mu^{r+1}:=\mu^1$. 
$\nu^i,\lambda^i$ ($1\leq i\leq r+1$) 
are defined similarly.
We define:
\begin{equation}\notag
l(\Vec{\mu})=\sum_{i=1}^{r}l(\mu^i),\qquad
\aut(\Vec{\mu})=\prod_{i=1}^r\aut({\mu^i}),\qquad
z_{\Vec{\mu}}=\prod_{i=1}^r z_{\mu^i},
\qquad
[\Vec{\mu}]=\prod_{i=1}^r[\mu^i].
\end{equation}
$l(\Vec{\nu}),l(\Vec{\lambda})$,
$\aut(\Vec{\nu}),\aut(\Vec{\lambda})$,
$z_{\Vec{\nu}},z_{\Vec{\lambda}}$ and
$[\Vec{\nu}],[\Vec{\lambda}]$
are defined in the same manner.

A triple $(\Vec{\mu},\Vec{\nu},\Vec{\lambda})$
is an {\em $r$-set}
if it satisfies
\begin{equation}\notag
|\mu^i|+|\lambda^i|=|\nu^i|+|\lambda^{i+1}|\qquad 
(1\leq {}^{\forall}i\leq r).
\end{equation}
$(|\mu^i|+|\lambda^{i}|)_{1\leq i\leq r}$ is called the degree
of the $r$-set.

As the consequence of 
the above lemma,
the partition function is 
written in terms of the matrix elements.
\begin{prop} \label{zop}
%
\begin{equation}\notag
{\mathcal Z}^{\gamma}_{\Vec{d}}(q)
=
(-1)^{{\gamma}\cdot\Vec{d}}
\sum_{\begin{subarray}{c}
(\Vec{\mu},\Vec{\nu},\Vec{\lambda});
\\
\text{$r$-set of}\\ 
\text{degree $\Vec{d}$}
\end{subarray}}
\frac{(-1)^{l(\Vec{\mu})+l(\Vec{\nu})}}
{z_{\Vec{\mu}}z_{\Vec{\nu}}z_{\Vec{\lambda}}
[\Vec{\mu}][\Vec{\nu}]}
\prod_{i=1}^r
\langle \lambda^i\cup \mu^i|
q^{(\gamma_i+2){\mathcal F}_2}|
\nu^i\cup\lambda^{i+1}\rangle.
\end{equation}
\end{prop}%
\begin{proof}
Using lemmas \ref{wmunu} and (\ref{hk}),
we write ${\mathcal Z}_{\Vec{d}}^{\gamma}(q)$ as follows.
\begin{equation}\notag
\begin{split}
(-1)^{\gamma\cdot\Vec{d}}{\mathcal Z}_{\Vec{d}}^{\gamma}(q)
=
\sum_{\begin{subarray}{c}
(\Vec{\mu},\Vec{\nu},\Vec{\eta});
\\
\text{$r$-set of}\\
\text{degree $\Vec{d}$}
\end{subarray}}
\frac{(-1)^{l(\Vec{\mu})+l(\Vec{\nu})}}
{z_{\Vec{\mu}}z_{\Vec{\nu}}z_{\Vec{\eta}}
[\Vec{\mu}][\Vec{\nu}]
}
\prod_{i=1}^r
\Big(
\sum_{\lambda^i;|\lambda^i|=d_i}
\langle \mu^i\cup\eta^i|q^{{\mathcal F}_2}|v_{\lambda^i}\rangle
\langle v_{\lambda^i}|q^{\gamma_i{\mathcal F}_2}|v_{\lambda^i}\rangle
\langle v_{\lambda^i}|q^{{\mathcal F}_2}|\nu^i\cup\eta^{i+1}\rangle
\Big).
\end{split}
\end{equation}
The first factor in the bracket comes from
$W_{\lambda^{i-1},\lambda^i}(q)$, 
the second is the factor $q^{\gamma_i{{\mathcal F}_2}}$
and 
the third comes from $W_{\lambda^i,\lambda^{i+1}}(q)$.
The bracket term is equal to
\begin{equation}\notag
\langle \mu^i\cup\eta^i|
q^{(\gamma_i+2){\mathcal F}_2}
|\nu^i\cup\eta^{i+1}\rangle
.
\end{equation}
If we replace the letter $\eta$ with $\lambda$,
we obtain the proposition.
\end{proof}
\section{Partition Function as Graph Amplitudes}\label{me}
In this section, we
express the partition function as the sum
of some amplitudes over possibly disconnected graphs.
(The term an {\em amplitude} will be used
to denote a map from a set of graphs to a ring.)

The graph description 
is achieved in two steps.
Firstly,
in subsections \ref{GD} and \ref{ME},
we study the matrix element 
appeared in the partition function:
\begin{equation}\label{me0}
\langle{\mu}|q^{a{\mathcal F}_2}|{\nu}\rangle
\qquad (a\in \mathbb{Z});
\end{equation}
we introduce the graph set 
associated to  $\mu,\nu$ and $a$
and describe the matrix element
as the sum of certain amplitude 
over the set.
Secondly, in subsection \ref{CA}
we turn to the whole partition function;
we combine
these graph sets 
and the amplitude to make
another type of a graph set 
and amplitude;
then we rewrite
the partition function in
terms of them.

\subsection{Notations}
Let us  briefly summarize  the notations on graphs (see \cite{gy}). 
\begin{itemize}
\item
A {\em graph} $G$ is a pair of the vertex set $V(G)$ 
and the edge set $E(G)$ such that
every edge has two vertices associated to it.
We only deal with graphs whose vertex sets and edge sets 
are finite.
\item
A {\em directed} graph is a graph whose edges has directions.
\item
A {\em label} of a graph $G$ is a
map from the vertex set $V(G)$ or from 
the edge set $E(G)$ to a set.
\item
An {\em isomorphism} between two graphs $G$ and $F$ is
a pair of bijections $V(G)\to V(F)$
and $E(G)\to E(F)$ preserving the incidence relations.
An {\em isomorphism} of labeled graphs is 
a graph isomorphism that preserves the labels.
\item
An {\em automorphism} is an isomorphism of the (labeled) 
graph to itself.
\item
The {\em graph union} of two graphs $G$ and $F$ is the graph 
whose vertex set and edge set are the disjoint unions,
respectively, of the vertex sets 
$V(G)$, $V(F)$ and
of the edge sets $E(G),E(F)$. 
It is denoted by $G\cup F$.
\item
If $G$ is connected,
the {\em cycle rank} (or {\em Betti number}) is $\#E(G)-\#V(G)+1$.
If $G$ is not connected,
its cycle rank is the sum of those of connected components.
\item
A connected graph with the cycle rank zero is 
a {\em tree} 
and the graph union of trees  is a {\em forest}.
\end{itemize}

We also use the following notations.
\begin{itemize}
%
\item
A set of not necessarily connected graphs is denoted by
a symbol with the superscript $\bullet$ and the subset
of connected graphs  by the same symbol with 
$\circ$ (e.g. $\mathbb{G}^{\bullet}$ and $\mathbb{G}^{\circ}$).
\item 
For any finite set of integers $s=(s_1,s_2,\ldots,s_l)$,
\begin{equation}\notag
|s|=\sum_i s_i.
\end{equation}
Note that $|s|$ can be negative if $s$ has negative elements. 

When $s$ has at least one nonzero element, we define
\begin{equation}\notag
{\rm gcd}(s)=\text{greatest common divisor of }\{|s_i|,s_i\neq 0\}
\end{equation}
where $|s_i|$ is the absolute value of $s_i$.
\end{itemize}

\subsection{Graph Description of VEV} \label{GD}
%
By the relation (\ref{f2alpha}), 
the matrix element (\ref{me0}) is rewritten as
the vacuum expectation value:
\begin{equation}\begin{split}\label{me1}
\langle{\mu}|q^{a{\mathcal F}_2}|{\nu}\rangle
&=
\langle{\mathcal E}_{\mu_{l(\mu)}}(0)\cdots
{\mathcal E}_{\mu_1}(0)
{\mathcal E}_{-\nu_1}(a\nu_1)\cdots
{\mathcal E}_{-\nu_{l(\nu)}}(a\nu_{l(\nu)})
\rangle.
\end{split}
\end{equation}
%
%
Therefore,
in this subsection, we consider
the  VEV
\begin{equation}\label{vev}
\langle {\mathcal E}_{c_1}(n_1)
{\mathcal E}_{c_2}(n_2)\cdots
{\mathcal E}_{c_l}(n_l)
\rangle
\end{equation}
where
\begin{equation}\notag
\Vec{c}=(c_1,\ldots,c_l),\qquad
\Vec{n}=(n_1,\ldots,n_l)
\end{equation}
is a pair of ordered sets of integers of the same length $l$.
We assume that 
$(c_i,n_i)\neq (0,0)$ ($1\leq{^\forall} i\leq l$)
because ${\mathcal E}_{0}(0)$ is not well-defined.
We also assume that
$|\Vec{c}|=0$
because (\ref{vev}) vanishes otherwise.

\subsubsection{Set of Graphs}
When we compute the VEV (\ref{vev}),
we 
make use of the commutation relation
several times.
We associate to this process
a graph generating algorithm in a natural way.
\begin{enumerate}%

\item[0.]
In the computation,
we start with (\ref{vev});
there are $l$ ${\mathcal E}_{c_i}(n_i)$'s in the angle bracket
and
we associate to each a black vertex;
we
draw $l$ black vertices horizontally on $\mathbb{R}^2$
and
assign $i$ and $(c_i,n_i)$ to the $i$-th vertex 
(counted from the left).
The picture is as follows.
\begin{center}
\unitlength .15cm
\begin{picture}(20,10)(-5,0)
\thicklines
\put(-2,3){\circle*{1}}
\put(5,3){\circle*{1}}
\multiput(10,3)(1,0){5}{\line(1,0){.1}}
\put(20,3){\circle*{1}}
\put(-2.5,7){1}
\put(-5,5){$(c_1,n_1)$}
\put(4.5,7){2}
\put(2.5,5){$(c_2,n_2)$}
\put(19.5,7){$l$}
\put(17.5,5){$(c_l,n_l)$}
\end{picture}
\end{center}

\item 
The 
first step in the computation is 
to take the rightmost adjacent pair
${\mathcal E}_{a}(m),{\mathcal E}_b(n)$
with $a\geq 0$ and $b< 0$ 
and  apply the commutation relation 
\begin{equation}\notag
{\mathcal E}_a(m){\mathcal E}_b(n)=
{\mathcal E}_b(n){\mathcal E}_a(m)
+
\begin{cases}
\Bigg[\,
\begin{vmatrix}a&m\\b&n\\\end{vmatrix}\,\Bigg]
{\mathcal E}_{a+b}(m+n)
& (a+b,m+n)\neq (0,0)\\
a&(a+b,m+n)=(0,0)
\end{cases}
\end{equation}
On the graph side,
we express this as follows.
%
\newcommand{\before}{  
\unitlength .15cm
\begin{picture}(7,8)(-3,-2)
\thicklines
\put(0,5){\circle*{1}}
\put(5,5){\circle*{1}}
\put(-3,7){$(a,m)$}
\put(3,7){$(b,n)$}
\end{picture}
}
\newcommand{\crossing}{
\unitlength .15cm
\begin{picture}(7,8)(-3,-2)
\thicklines
\put(0,5){\line(1,-1){5}}
\put(5,5){\line(-1,-1){5}}
\put(0,5){\circle*{1}}
\put(5,5){\circle*{1}}
\put(-3,7){$(a,m)$}
\put(3,7){$(b,n)$}
\put(0,0){\circle*{1}}
\put(5,0){\circle*{1}}
\put(-3,-2){$(b,n)$}
\put(3,-2){$(a,m)$}
\end{picture}
}
\newcommand{\collide}{
\unitlength .15cm
\begin{picture}(7,8)(-3,-3)
\thicklines
\put(0,5){\line(1,-2){2.5}}
\put(5,5){\line(-1,-2){2.5}}
\put(0,5){\circle*{1}}
\put(5,5){\circle*{1}}
\put(-3,7){$(a,m)$}
\put(3,7){$(b,n)$}
\put(2.5,0){\circle*{1}}
\put(-2,-2){$(a+b,m+n)\neq (0,0)$}
\end{picture}
}
\newcommand{\coll}{
\unitlength .15cm
\begin{picture}(7,10)(-3,-2)
\thicklines
\put(0,5){\line(1,-2){2.5}}
\put(5,5){\line(-1,-2){2.5}}
\put(0,5){\circle*{1}}
\put(5,5){\circle*{1}}
\put(-3,7){$(a,m)$}
\put(3,7){$(b,n)$}
\put(2.5,0){\circle{1}}
\put(-1,-2){$(a+b,m+n)=(0,0)$}
\end{picture}
}
\begin{equation}\notag
\raisebox{-.7cm}{\before} 
\qquad = \qquad
\raisebox{-.7cm}{\crossing} 
\qquad + \qquad
\begin{cases}
\raisebox{-.7cm}{\collide}\\
\\
\raisebox{-.5cm}{\coll}
\end{cases}
\end{equation}
For other black vertices, we do as follows.
\newcommand{\oth}{
\unitlength .15cm
\begin{picture}(5,8)(-2,-2)
\thicklines
\put(0,5){\circle*{1}}
\put(-2,7){$(c,h)$}
\end{picture}
}
\newcommand{\othe}{
\unitlength .15cm
\begin{picture}(10,8)(-4,-2)
\thicklines
\put(0,5){\line(0,-1){5}}
\put(0,5){\circle*{1}}
\put(0,0){\circle*{1}}
\put(-3,7){$(c,h)$}
\put(-3,-3){$(c,h)$}
\end{picture}
}
\begin{equation}\notag
\raisebox{-.6cm}{\oth} \qquad\qquad \Longrightarrow 
\qquad\qquad \raisebox{-.5cm}{\othe}
\end{equation}

\item
Secondly,
we use the following relations if applicable:
\begin{equation}\notag
\begin{split}
&\langle\cdots{\mathcal E}_a(m)\rangle
=
\begin{cases}0&(a>0)\\ 
\langle\cdots\rangle \frac{1}{[m]}&(a=0)
\end{cases}
\\
&\langle{\mathcal E}_b(n)\cdots\rangle=0\quad (b<0),
\qquad \text{ and }\langle 1\rangle=1.
\end{split}
\end{equation}
Accordingly, in the drawing, we do as follows.
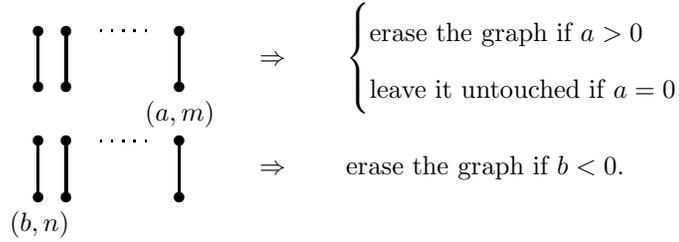
\newcommand{\hor}{
\unitlength .15cm
\begin{picture}(21,8)(0,-2)
\thicklines
\put(10,5){\line(0,-1){5}}
\put(7.5,5){\line(0,-1){5}}
\put(20,5){\line(0,-1){5}}
%
\put(7.5,5){\circle*{1}}
\put(10,5){\circle*{1}}
\multiput(13,5)(1,0){5}{\line(1,0){.1}}
\put(20,5){\circle*{1}}
\put(7.5,0){\circle*{1}}
\put(10,0){\circle*{1}}
\put(20,0){\circle*{1}}
\put(17,-3){$(a,m)$}
\end{picture}
}
\newcommand{\hori}{
\unitlength .15cm
\begin{picture}(21,8)(0,-2)
\thicklines
\put(10,5){\line(0,-1){5}}
\put(7.5,5){\line(0,-1){5}}
\put(20,5){\line(0,-1){5}}
\put(7.5,5){\circle*{1}}
\put(10,5){\circle*{1}}
\multiput(13,5)(1,0){5}{\line(1,0){.1}}
\put(7.5,0){\circle*{1}}
\put(10,0){\circle*{1}}
\put(20,0){\circle*{1}}
\put(20,5){\circle*{1}}
\put(5,-3){$(b,n)$}
\end{picture}
}
\begin{equation}\notag
\begin{split}
&\raisebox{-.6cm}{\hor} \qquad \Rightarrow
\qquad
\begin{cases}
\text{erase the graph if $a>0$}\\
\text{leave it untouched if $a=0$}
\end{cases}
\\
&\raisebox{-.6cm}{\hori} \qquad \Rightarrow
\qquad
\text{erase the graph if $b<0$}.
\end{split}
\end{equation} 
\item
We
repeat these steps 
until all terms become constants.
Although the number of terms may increase at first, 
it  stabilizes in the end and the computation stops with finite 
processes.
\item
Finally, 
we simplify the drawings: 
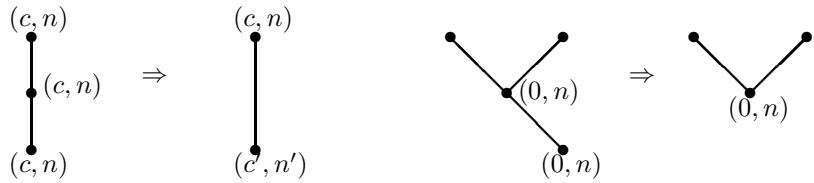
\newcommand{\extravertex}{
\unitlength .15cm
\begin{picture}(10,13)(-4,-4)
\thicklines
\put(0,0){\line(0,1){10}}
\put(0,0){\circle*{1}}
\put(0,5){\circle*{1}}
\put(0,10){\circle*{1}}
\put(-2,11){$(c,n)$}
\put(1,5){$(c,n)$}
\put(-2,-2){$(c,n)$}
\end{picture}
}
\newcommand{\deletevertex}{
\unitlength .15cm
\begin{picture}(10,13)(-4,-4)
\thicklines
\put(0,0){\line(0,1){10}}
\put(0,0){\circle*{1}}
\put(0,10){\circle*{1}}
\put(-2,11){$(c,n)$}
\put(-2,-2){$(c',n')$}
\end{picture}
}
\newcommand{\extraroot}{
\unitlength .15cm
\begin{picture}(10,13)(-15,-4)
\thicklines
\put(0,0){\line(-1,1){10}}
\put(-5,5){\line(1,1){5}}
\put(0,0){\circle*{1}}
\put(-5,5){\circle*{1}}
\put(-10,10){\circle*{1}}
\put(0,10){\circle*{1}}
\put(-4,4.5){$(0,n)$}
\put(-2,-2){$(0,n)$}
\end{picture}
}
\newcommand{\deleteroot}{
\unitlength .15cm
\begin{picture}(5,13)(-12,-4)
\thicklines
\put(-5,5){\line(-1,1){5}}
\put(-5,5){\line(1,1){5}}
\put(-5,5){\circle*{1}}
\put(-10,10){\circle*{1}}
\put(0,10){\circle*{1}}
\put(-7,3){$(0,n)$}
\end{picture}
}

\begin{equation}\notag
\!\!\!\!\!\!\!\!\!\!\!\!\!\!\!
\raisebox{-1.5cm}{\extravertex}
               \quad \Rightarrow               
\quad \raisebox{-1.5cm}
{\deletevertex}
\qquad 
\raisebox{-1.5cm}{\extraroot}
  \qquad\qquad \Rightarrow               
\raisebox{-1.5cm}{\deleteroot}
\end{equation}
\end{enumerate}%

For concreteness,
we show the case of 
\begin{equation}\notag
\Vec{c}=(c,c,-c,-c),\qquad
\Vec{n}=(0,0,0,d)\qquad
(c>0,d\neq 0).
\end{equation}
We will omit the labels and irrelevant middle vertices.
We start with
\newcommand{\stone}{\raisebox{-.4cm}{
\unitlength .15cm
\begin{picture}(15,5)(-5,0)
\thicklines
\put(0,3){\circle*{1}}
\put(5,3){\circle*{1}}
\put(10,3){\circle*{1}}
\put(15,3){\circle*{1}}
\end{picture}
}}
\begin{equation}\notag
\langle{\mathcal E}_{c}(0)
\underbrace{
{\mathcal E}_{c}(0){\mathcal E}_{-c}(0)
}
{\mathcal E}_{-c}(d)
\rangle
\qquad 
\stone
\end{equation}
We pick the pair with the underbrace and 
apply the step 1:
\newcommand{\sttwo}{\raisebox{-.3cm}{
\unitlength .15cm
\begin{picture}(15,6)(0,0)
\thicklines
\put(0,0){\line(0,1){5}}
\put(15,0){\line(0,1){5}}
\put(5,0){\line(1,1){5}}
\put(10,0){\line(-1,1){5}}
\put(0,0){\circle*{1}}
\put(5,0){\circle*{1}}
\put(10,0){\circle*{1}}
\put(15,0){\circle*{1}}
\put(0,5){\circle*{1}}
\put(5,5){\circle*{1}}
\put(10,5){\circle*{1}}
\put(15,5){\circle*{1}}
\end{picture}
}}
\newcommand{\stthree}{\raisebox{-.3cm}{
\unitlength .15cm
\begin{picture}(15,6)(0,0)
\thicklines
\put(0,0){\line(0,1){5}}
\put(15,0){\line(0,1){5}}
\put(5,5){\line(1,-2){2.5}}
\put(10,5){\line(-1,-2){2.5}}
\put(0,0){\circle*{1}}
\put(15,0){\circle*{1}}
\put(0,5){\circle*{1}}
\put(5,5){\circle*{1}}
\put(10,5){\circle*{1}}
\put(15,5){\circle*{1}}
\put(7.5,0){\circle{1}}
\end{picture}
}}
\begin{equation}\notag
\begin{split}
&\langle
{\mathcal E}_{c}(0){\mathcal E}_{-c}(0)
\underbrace{
{\mathcal E}_{c}(0){\mathcal E}_{-c}(d)
}
\rangle 
+
c\langle
\underbrace{
{\mathcal E}_{c}(0){\mathcal E}_{-c}(d)
}
\rangle
\\&
\sttwo
\qquad\qquad \qquad 
\stthree
\end{split}
\end{equation}
We skip the step 2 
since there is no place to apply the relation.  
We proceed to the step 1 again
and
take the commutation relation of the underbraced
pairs:
\newcommand{\stfour}{
\unitlength .15cm
\begin{picture}(15,7)(0,-1)
\thicklines
\put(0,5){\line(0,-1){5}}
\put(15,5){\line(-1,-1){5}}
\put(5,5){\line(2,-1){10}}
\put(10,5){\line(-1,-1){5}}
\put(0,0){\circle*{1}}
\put(5,0){\circle*{1}}
\put(10,0){\circle*{1}}
\put(15,0){\circle*{1}}
\put(0,5){\circle*{1}}
\put(5,5){\circle*{1}}
\put(10,5){\circle*{1}}
\put(15,5){\circle*{1}}
\end{picture}}
\newcommand{\stfive}{
\unitlength .15cm
\begin{picture}(15,7)(0,-1)
\thicklines
\put(0,5){\line(0,-1){5}}
\put(15,5){\line(-1,-1){5}}
\put(5,5){\line(1,-1){5}}
\put(10,5){\line(-1,-1){5}}
\put(0,0){\circle*{1}}
\put(5,0){\circle*{1}}
\put(10,0){\circle*{1}}
\put(0,5){\circle*{1}}
\put(5,5){\circle*{1}}
\put(10,5){\circle*{1}}
\put(15,5){\circle*{1}}
\end{picture}
}
\newcommand{\stsix}{
\unitlength .15cm
\begin{picture}(15,15)(0,4)
\thicklines
\put(0,15){\line(1,-1){10}}
\put(15,15){\line(-1,-1){10}}
\put(5,15){\line(1,-1){2.5}}
\put(10,15){\line(-1,-1){2.5}}
\put(0,15){\circle*{1}}
\put(5,15){\circle*{1}}
\put(10,15){\circle*{1}}
\put(15,15){\circle*{1}}
\put(10,5){\circle*{1}}
\put(5,5){\circle*{1}}
\put(7.5,12.5){\circle{1}}
\end{picture}
}
\newcommand{\stseven}{
\unitlength .15cm
\begin{picture}(15,15)(0,4)
\thicklines
\put(0,15){\line(1,-1){7.5}}
\put(15,15){\line(-1,-1){7.5}}
\put(5,15){\line(1,-1){2.5}}
\put(10,15){\line(-1,-1){2.5}}
\put(0,15){\circle*{1}}
\put(5,15){\circle*{1}}
\put(10,15){\circle*{1}}
\put(15,15){\circle*{1}}
\put(7.5,12.5){\circle{1}}
\put(7.5,7.5){\circle*{1}}
\end{picture}
}
\begin{equation}\notag
\begin{split}
&\langle
{\mathcal E}_{c}(0){\mathcal E}_{-c}(0)
{\mathcal E}_{-c}(d){\mathcal E}_{c}(0)
\rangle
+
[cd]\langle
{\mathcal E}_c(0)
{\mathcal E}_{-c}(0)
{\mathcal E}_{0}(d)
\rangle
+c\langle{\mathcal E}_{-c}(d){\mathcal E}_c(0)\rangle
+
c[cd]\langle{\mathcal E}_0(d)\rangle
\\
&\raisebox{.5cm}{\stfour}
\qquad \qquad \raisebox{.5cm}{\stfive}
\qquad \qquad{\stsix}\qquad\qquad 
\stseven
\end{split}
\end{equation}
This time we apply the step 2.
The first and the third term become zero.
In the second and the fourth terms,
${\mathcal E}_0(d)$'s are
replaced by $1/[d]$.
The result is
\begin{equation}\notag
\begin{split}
&
\frac{[cd]}{[d]}\langle
\underbrace{
{\mathcal E}_{c}(0){\mathcal E}_{-c}(0)
}
\rangle 
\qquad
+\qquad
\frac{c[cd]}{[d]}
\\&
\raisebox{.5 cm}{\stfive} \qquad \qquad \stseven
\end{split}
\end{equation}
Applying the step 1 again,
we obtain
\newcommand{\steight}{
\unitlength .15cm
\begin{picture}(15,7)(0,-1)
\thicklines
\put(0,5){\line(1,-1){5}}
\put(15,5){\line(-1,-1){5}}
\put(5,5){\line(1,-1){5}}
\put(10,5){\line(-2,-1){10}}
\put(0,0){\circle*{1}}
\put(5,0){\circle*{1}}
\put(10,0){\circle*{1}}
\put(0,5){\circle*{1}}
\put(5,5){\circle*{1}}
\put(10,5){\circle*{1}}
\put(15,5){\circle*{1}}
\end{picture}
}
\newcommand{\stnine}{
\unitlength .15cm
\begin{picture}(15,7)(0,-1)
\thicklines
\put(0,5){\line(1,-1){5}}
\put(15,5){\line(-1,-1){5}}
\put(5,5){\line(1,-1){5}}
\put(10,5){\line(-1,-1){5}}
\put(0,5){\circle*{1}}
\put(10,0){\circle*{1}}
\put(5,5){\circle*{1}}
\put(10,5){\circle*{1}}
\put(15,5){\circle*{1}}
\put(5,0){\circle{1}}
\end{picture}
}
\begin{equation}\notag
\begin{split}
&\frac{[cd]}{[d]}\langle
{\mathcal E}_{-c}(0){\mathcal E}_{c}(0)\rangle
\qquad+\qquad
\frac{c[cd]}{[d]}
\qquad+\qquad
\frac{c[cd]}{[d]}
\\
&\raisebox{.5cm}{\steight}\qquad\qquad
\raisebox{.5cm}{\stnine}\qquad\qquad
{\stseven}
\end{split}
\end{equation}
The first term is zero.
Hence all the terms become constants and 
the process is finished. 
What 
we obtained are the two nonzero terms and 
the corresponding two graphs:
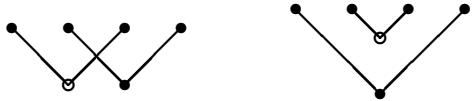
\begin{equation}\notag
\begin{split}
\raisebox{.5cm}{\stnine} \qquad\qquad \stseven
\end{split}
\end{equation}
%

\begin{definition} \label{defgcn}
${\rm Graph}^{\bullet}(\Vec{c}, \Vec{n})$
is  the set of graphs generated by the above 
recursion algorithm.
${\rm Graph}^{\circ}(\Vec{c}, \Vec{n})$
is the subset consisting of all connected graphs. 
\end{definition}%

Since
$F\in {\rm Graph}^{\bullet}(\Vec{c}, \Vec{n})$
is the graph union of trees,
we call $F$ a {\em VEV forest}.
A univalent vertex at the highest level is called 
{\em a leaf}.
The unique bivalent vertex at the lowest level of 
each connected component is 
called the {\em root}.

A VEV forest has two types of labels.
The labels attached to leaves are referred to as
the {\em leaf indices}.
The two component labels on every vertex $v$ is called the 
vertex-label of $v$ and denoted by $(c_v,n_v)$.

We 
define $\overline{F}$ to be the graph obtained from 
a VEV forest $F$ by forgetting
the leaf indices.
Two  VEV forests $F$ and $F'$ are {\em equivalent}
if $\overline{F}\cong\overline{F}'$.
The set of equivalence classes in
${\rm Graph}^{\bullet}(\Vec{c},\Vec{n})$ 
and
${\rm Graph}^{\circ}(\Vec{c},\Vec{n})$)
are denoted by
$\overline{{\rm Graph}}^{\bullet}(\Vec{c},\Vec{n})$
and
$\overline{{\rm Graph}}^{\circ}(\Vec{c},\Vec{n})$.
A connected component $T$ of $\overline{F}$
is called a {\em VEV tree}.
Note that $T$ is regarded as an element
of
$\overline{{\rm Graph}}^{\circ}({\Vec{c}}',\Vec{n}')$
with some $({\Vec{c}}'\Vec{n}')$.

Let $V_2(F)$ be the set of vertices
which have two adjacent vertices
at the upper level.
For a connected component $T$ of $F$, $V_2(T)$ is defined similarly.
$L(v)$ and $R(v)$ denote the upper left and right
vertices adjacent to $v\in V_2(F)$.
\newcommand{\LR}{
\unitlength .15cm
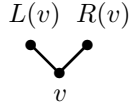
\begin{picture}(10,8)(0,0)
\thicklines
\put(5,5){\line(1,-1){2.5}}
\put(10,5){\line(-1,-1){2.5}}
\put(3,7){$L(v)$}
\put(5,5){\circle*{1}}
\put(10,5){\circle*{1}}
\put(9,7){$R(v)$}

\put(7.5,2.5){\circle*{1}}
\put(7,0){$v$}
\end{picture}
}
\begin{equation}\notag
\LR
\end{equation}
A left leaf (right leaf) is a
leaf which is  $L(v)$ ($R(v)$) of some vertex $v$.

We summarize the properties of 
the vertex labels of a VEV tree $T$.
\begin{enumerate}
\item 
$c_{\text{root}}=
\sum_{v:\text{leaves}}c_{v}=0$.
\item If a vertex $v$ is white,
$(c_v,n_v)=(0,0)$ and it is the root vertex. 
\item
If $v$ is black, $(c_v,n_v)\neq (0,0)$. 
%
%
\item 
$c_{L(v)}\geq 0$ and $c_{R(v)}<0$ for $v\in V_2(T)$.
\item 
$c_{v}=c_{L(v)}+c_{R(v)}$ and $n_{v}=n_{L(v)}+n_{R(v)}$
for $v\in V_2(T)$.
\end{enumerate}

\subsubsection{Amplitude}
%
%
In the computation process,
the birth of 
each black vertex $v\in V_2(F)$ 
means that the corresponding term is multiplied by
the factor
\begin{equation}\label{zetav}
\begin{split}
\zeta_v=
\begin{vmatrix}c_{L(v)} &n_{L(v)}
\\c_{R(v)}&n_{R(v)}
\\\end{vmatrix}
&\qquad (v\in V_2(F)).
\end{split}
\end{equation}
Moreover if the vertex is the root, 
then the corresponding term is multiplied by $1/[n_v]$.
On the other hand, the birth of a white vertex $v$ 
means that the corresponding term is multiplied by
the factor $c_{L(v)}$.
Therefore 
we define the amplitude ${\mathcal A}(F)$
for $F\in {\rm Graph}^{\bullet}(\Vec{c}, \Vec{n})$
by
%
%
\begin{equation}\label{amp}
\begin{split}
{\mathcal A}(F)&=
\prod_{T:\text{ VEV tree in $F$}}
{\mathcal A}(T),
\\
{\mathcal A}(T)&=
\begin{cases}
{\prod_{v\in V_2(T)}[\zeta_v]}/
{[n_{\text{root}}]}
&(n_{\text{root}}\neq 0)
\\
c_{L(\text{root})}
\prod_{\begin{subarray}{c}v\in V_2(T),\\
     v\neq \text{root}\end{subarray}}
[\zeta_v]
&(n_{\text{root}}=0)
\end{cases}.
\end{split}
\end{equation}
From the definition of
${\rm Graph}^{\bullet}(\Vec{c}, \Vec{n})$,
it is clear that
the sum of all the amplitude ${\mathcal A}(F)$
is equal to
the VEV (\ref{vev}).
%
\begin{prop}%
\label{eprop}
\begin{equation}\notag
\langle {\mathcal E}_{c_1}(n_1)
{\mathcal E}_{c_2}(n_2)\cdots
{\mathcal E}_{c_l}(n_l)
\rangle
=\sum_{F\in {\rm Graph}^{\bullet}(\Vec{c},\Vec{ n})}
{\mathcal A}(F).
\end{equation}
\end{prop}
%

\subsection{Matrix Elements}\label{ME}
In this subsection, we apply the result of preceding subsection 
to the matrix element (\ref{me0}).

\subsubsection{Graphs}
Let $\mu,\nu$ be two partitions $|\mu|=|\nu|=d>1$, $a\in\mathbb{Z}$.
\begin{definition}%
${\rm Graph}_a^{\bullet}({\mu,\nu})$ is the set 
${\rm Graph}^{\bullet}(\Vec{c},\Vec{n})$ 
with
\begin{equation}\notag
\Vec{c}=
(\mu_{l(\mu)},\ldots,\mu_{1},
-\nu_{1},\ldots,-\nu_{l(\nu)}),\qquad
\Vec{n}=
(\underbrace{0,\ldots,0}_{l(\mu) \text{ times}},
a\nu_1,\ldots,a\nu_{l(\nu)}).
\end{equation}
The subset of connected graphs is denoted by
${\rm Graph}_a^{\circ}({\mu,\nu})$.
The set of equivalence classes 
of ${\rm Graph}_a^{\bullet}({\mu,\nu})$
and ${\rm Graph}_a^{\circ}({\mu,\nu})$
obtained by forgetting the leaf indices
are
denoted by 
$\overline{{\rm Graph}}_a^{\bullet}({\mu,\nu})$
and $\overline{{\rm Graph}}_a^{\circ}({\mu,\nu})$.
%
\end{definition}%

The next proposition
follows immediately from (\ref{me1}) and proposition \ref{eprop}.
\begin{prop}\label{me3}
\begin{equation}\notag
\langle{\mu}|q^{a{\mathcal F}_2}|{\nu}\rangle
=\sum_{F\in {\rm Graph}_a^{\bullet}(\mu,\nu)}
{\mathcal A}(F).
\end{equation}
\end{prop}
%
\begin{example}\label{exmunu}
Some examples of ${\rm Graph}_a^{\bullet}(\mu,\nu)$
and the amplitudes of its elements are shown below.

$\mu=(\mu_1,\mu_2,\mu_3)$, $\nu=(d)$ where $d=\mu_1+\mu_2+\mu_3$:
\newcommand{\muda}{
\unitlength .15cm
\begin{picture}(15,15)(0,4)
\thicklines
\put(0,15){\line(1,-1){7.5}}
\put(15,15){\line(-1,-1){7.5}}
\put(5,15){\line(1,-1){5}}
\put(10,15){\line(1,-1){2.5}}
\put(0,15){\circle*{1}}
\put(5,15){\circle*{1}}
\put(10,15){\circle*{1}}
\put(15,15){\circle*{1}}
\put(10,10){\circle*{1}}
\put(12.5,12.5){\circle*{1}}
\put(7.5,7.5){\circle*{1}}
\end{picture}
}
\newcommand{\mudzero}{
\unitlength .15cm
\begin{picture}(15,15)(0,4)
\thicklines
\put(0,15){\line(1,-1){7.5}}
\put(15,15){\line(-1,-1){7.5}}
\put(5,15){\line(1,-1){5}}
\put(10,15){\line(1,-1){2.5}}
\put(0,15){\circle*{1}}
\put(5,15){\circle*{1}}
\put(10,15){\circle*{1}}
\put(15,15){\circle*{1}}
\put(10,10){\circle*{1}}
\put(12.5,12.5){\circle*{1}}
\put(7.5,7.5){\circle{1}}
\end{picture}
}
\begin{equation}\notag
\begin{split}
a\neq 0:\qquad&\raisebox{-.8cm}{\muda}
\qquad\qquad \frac{[ad\mu_1][ad\mu_2][ad\mu_3]}{[ad]}
\\
a=0:\qquad&
\raisebox{-.8cm}{\mudzero}
\qquad\qquad 0
\end{split}
\end{equation}
$\mu=\nu=(c,c)$, $a\neq 0$:
\newcommand{\cca}{
\unitlength .15cm
\begin{picture}(15,5)(0,0)
\thicklines
\put(0,5){\line(1,-1){5}}
\put(15,5){\line(-1,-1){5}}
\put(5,5){\line(1,-1){5}}
\put(10,5){\line(-1,-1){5}}
\put(0,5){\circle*{1}}
\put(5,5){\circle*{1}}
\put(10,5){\circle*{1}}
\put(15,5){\circle*{1}}
\put(5,0){\circle*{1}}
\put(10,0){\circle*{1}}
\end{picture}
}
\newcommand{\ccatwo}{
\unitlength .15cm
\begin{picture}(15,15)(0,7.5)
\thicklines
\put(0,15){\line(1,-1){7.5}}
\put(15,15){\line(-1,-1){7.5}}
\put(5,15){\line(1,-1){2.5}}
\put(10,15){\line(-1,-1){2.5}}
\put(0,15){\circle*{1}}
\put(5,15){\circle*{1}}
\put(10,15){\circle*{1}}
\put(15,15){\circle*{1}}
\put(7.5,12.5){\circle*{1}}
\put(7.5,7.5){\circle*{1}}
\end{picture}
}
\newcommand{\ccasan}{
\unitlength .15cm
\begin{picture}(15,15)(0,7.5)
\thicklines
\put(0,15){\line(1,-1){7.5}}
\put(15,15){\line(-1,-1){7.5}}
\put(5,15){\line(1,-1){5}}
\put(10,15){\line(-1,-1){2.5}}
\put(0,15){\circle*{1}}
\put(5,15){\circle*{1}}
\put(10,15){\circle*{1}}
\put(15,15){\circle*{1}}
\put(10,10){\circle*{1}}
\put(7.5,12.5){\circle*{1}}
\put(7.5,7.5){\circle*{1}}
\end{picture}
}
\begin{equation}\notag
\begin{split}
&\ccasan \qquad \raisebox{.4cm}{\cca} \qquad \ccatwo
\\
&\quad\frac{[ac^2]^2[2ac^2]}{[ac]}
\qquad\qquad\quad
\Big(\frac{[ac^2]}{[ac]}\Big)^2
\qquad\qquad
\Big(\frac{[ac^2]}{[ac]}\Big)^2
\end{split}
\end{equation}
$\mu=\nu=(c,c)$ and $a=0$.
\newcommand{\cczero}{
\unitlength .15cm
\begin{picture}(15,0)(0,0)
\thicklines
\put(0,5){\line(1,-1){5}}
\put(15,5){\line(-1,-1){5}}
\put(5,5){\line(1,-1){5}}
\put(10,5){\line(-1,-1){5}}
\put(0,5){\circle*{1}}
\put(5,5){\circle*{1}}
\put(10,5){\circle*{1}}
\put(15,5){\circle*{1}}
\put(5,0){\circle{1}}
\put(10,0){\circle{1}}
\end{picture}
}
\newcommand{\cczerotwo}{
\unitlength .15cm
\begin{picture}(15,15)(0,7.5)
\thicklines
\put(0,15){\line(1,-1){7.5}}
\put(15,15){\line(-1,-1){7.5}}
\put(5,15){\line(1,-1){2.5}}
\put(10,15){\line(-1,-1){2.5}}
\put(0,15){\circle*{1}}
\put(5,15){\circle*{1}}
\put(10,15){\circle*{1}}
\put(15,15){\circle*{1}}
\put(7.5,12.5){\circle{1}}
\put(7.5,7.5){\circle{1}}
\end{picture}
}
\begin{equation}\notag
\begin{split}
&\raisebox{.2cm}{\cczero} \qquad \cczerotwo
\\
&\qquad\quad c^2
\qquad\qquad\qquad\qquad
c^2
\end{split}
\end{equation}

\end{example}

\subsection{Combined Amplitude}\label{CA}
%
In the previous section, 
we have written the matrix elements as 
the sum of  amplitudes.
In section \ref{pfof},
we have expressed
the partition function 
${\mathcal Z}_{\Vec{d}}^{\gamma}(q)$ 
in terms of matrix elements.
Combining these results, we obtain
\begin{equation}\notag
{\mathcal Z}^{\gamma}_{\Vec{d}}(q)
=
(-1)^{{\gamma}\cdot\Vec{d}}
\sum_{\begin{subarray}{c}
(\Vec{\mu},\Vec{\nu},\Vec{\lambda});\\
\text{$r$-set of}\\
\text{degree $\Vec{d}$}
\end{subarray}}
\frac{(-1)^{l(\Vec{\mu})+l(\Vec{\nu})}}
{z_{\Vec{\lambda}}z_{\Vec{\mu}}z_{\Vec{\nu}}[\Vec{\mu}][\Vec{\nu}]}
\sum_{\begin{subarray}{c}(F_i)_{1\leq i\leq r};\\
F_i\in {\rm Graph}^{\bullet}_{(\gamma_i+2)}
(\mu^i\cup\lambda^i,\nu^i\cup\lambda^{i+1})
\end{subarray}
}
\prod_{i=1}^r
{\mathcal A}(F_i).
\end{equation}
Next,
we will associate the each term in the RHS 
a new type of graph and
its amplitudes.

\subsubsection{Combined Forests}
Let $(\Vec{\mu},\Vec{\nu},\Vec{\lambda})$ be 
an $r$-set.
We construct a new graph from an
$r$-tuple 
\begin{equation}\notag
(F_1,\ldots,F_r),
\qquad F_i\in 
{\rm Graph}^{\bullet}_{(\gamma_i+2)}
(\mu^i\cup\lambda^i,\nu^i\cup\lambda^{i+1})
\end{equation}
as follows.
\begin{enumerate}
\item
Assign the label $i$ to each $F_i$ and make the graph union.
\item
Recall that
every left leaf of $F_i$ is associated to
a part of $\lambda^i$ or $\mu^i$ 
and every right leaf to a part of $\nu^i$ or $\lambda^{i+1}$
(see the remark below).
Since there are two leaves 
associated to $\lambda^i_j$ $(1\leq j\leq l(\lambda^i))$,
the one in $F_{i-1}$ and the other in $F_{i}$,
join them and assign the label $\lambda^{i}_j$ to the new edge.
\end{enumerate}
The resulting graph $W$ is a set of VEV forests marked by $1,\ldots,r$
and
joined through leaves.
We call $W$ a {\em combined forest}.
The new edges are called the {\em bridges}.
The label of a bridge $b$ is denoted by $h(b)$.
Note that the number of bridges in $W$ 
is equal to $l(\Vec{\lambda})$.

\begin{remark}%
To be precise,
every left 
leaf corresponds to 
a part of $\mu^i\cup\lambda^i$.
So
a problem would arise
if there are equal parts in $\mu^i$ and $\lambda^i$
since they are indistinguishable in $\mu^i\cup\lambda^i$.
To solve it, we set the following rule.
Assume $\mu^i$ and $\lambda^i$ 
have an equal part, say $k$.
Then
there are several left  leaves with
the vertex label $(k,0)$.
We regard the outer such leaves
come from $\lambda^{i}$ 
and the inner leaves
come from $\mu^i$. 
For right leaves, we set the same rule;
outer leaves
come from $\lambda^{i+1}$
and the inner leaves
come from $\nu^i$.
\end{remark}
\begin{definition}%
The set ${\rm Comb}^{\bullet}_{\gamma}
(\Vec{\mu},\Vec{\nu},\Vec{\lambda})$
is 
the set of combined forests constructed by the above procedure.
The set ${\rm Comb}_{\gamma}^{\circ}
(\Vec{\mu},\Vec{\nu},\Vec{\lambda})$
is 
the subset consisting of all connected combined forests. 
\end{definition}%

\begin{figure}[ht]
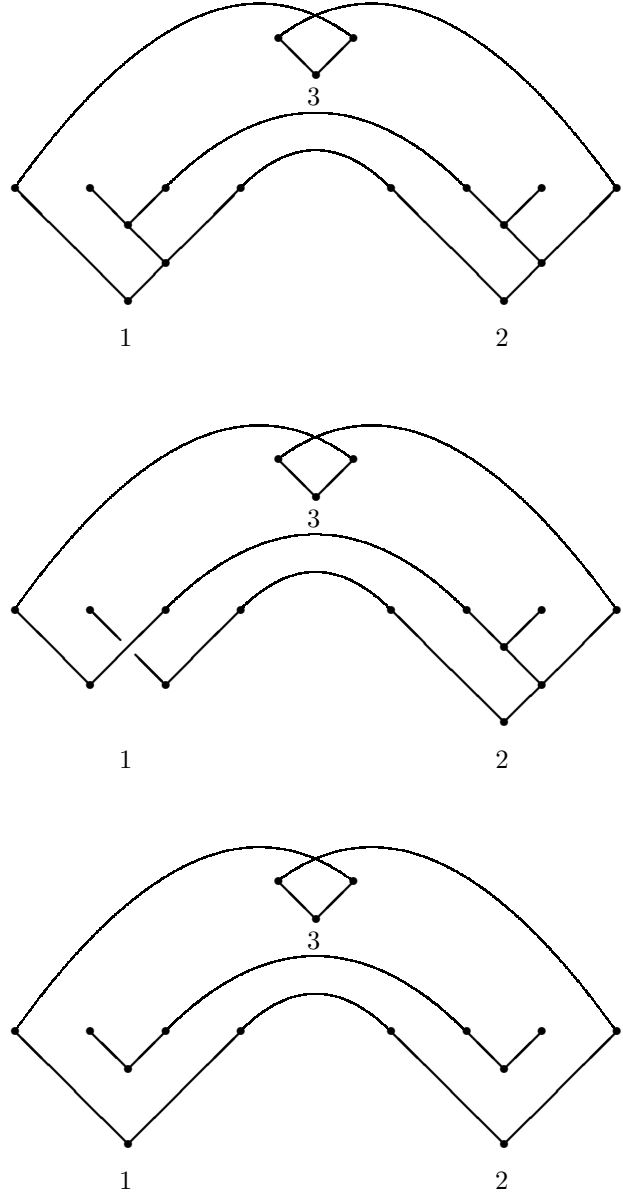

\begin{equation}\notag
\!\!\!\!\!\!\!\!\!\!
\!\!\!\!\!\!\!\!\!\!
\!\!\!\!\!\!\!\!\!\!
\!\!\!\!\!\!\!\!\!\!
\!\!\!\!\!\!\!\!\!\!
\!\!\!\!\!\!\!\!\!\!
\!\!\!\!\!\!\!\!\!\!
\!\!\!\!\!\!\!\!\!\!
\comb{\treea}{\treea}
\end{equation}
\begin{equation}\notag
\!\!\!\!\!\!\!\!\!\!
\!\!\!\!\!\!\!\!\!\!
\!\!\!\!\!\!\!\!\!\!
\!\!\!\!\!\!\!\!\!\!
\!\!\!\!\!\!\!\!\!\!
\!\!\!\!\!\!\!\!\!\!
\!\!\!\!\!\!\!\!\!\!
\!\!\!\!\!\!\!\!\!\!
\comb{\treeb}{\treea}
\end{equation}
\begin{equation}\notag
\!\!\!\!\!\!\!\!\!\!
\!\!\!\!\!\!\!\!\!\!
\!\!\!\!\!\!\!\!\!\!
\!\!\!\!\!\!\!\!\!\!
\!\!\!\!\!\!\!\!\!\!
\!\!\!\!\!\!\!\!\!\!
\!\!\!\!\!\!\!\!\!\!
\!\!\!\!\!\!\!\!\!\!
\comb{\treec}{\treec}
\end{equation}
\caption{Examples of combined forests.
The graphs shown are  three of nine elements
in the set ${\rm Comb}_{\gamma}^{\bullet}
(\vec{\mu},\vec{\nu},\vec{\lambda})$
where
$r=3$, $\gamma_i+2\neq 0$ and
$\lambda^1=(1)=\lambda^3,\lambda^2=(1,1)$,
$\mu^1=(1)=\nu^2, \mu^2=\mu^3=\nu^1=\nu^3=\emptyset$.
The leaf indices and 
the vertex labels are omitted.
}
\label{excomb1}
\end{figure}

We show some examples in the figure \ref{excomb1}.
%

\subsection{Combined Amplitude and Partition Function}

We first define a slightly different 
amplitude ${\mathcal B}$.
For a VEV tree
$T\in \overline{{\rm Graph}}_{a}^{\circ}(\mu,\nu)$,
\begin{equation}\notag
{\mathcal B}(T)=\frac{{\mathcal A}(T)}{[\mu][\nu]}.
\end{equation}
We
define the amplitude ${\mathcal H}$ of 
a combined forest
$W
\in {\rm Comb}_{\gamma}^{\bullet}
(\vec{\mu},\vec{\nu},\vec{\lambda})$
by
\begin{equation}\notag
{\mathcal H}(W)=
(-1)^{L_1(W)+L_2(W)}
\prod_{T:\text{ VEV tree in $W$}}{\mathcal B}(T)
\prod_{b;\text{bridge}}[h(b)]^2
\end{equation}
where $L_1(W)=l(\Vec{\mu})+l(\Vec{\nu})$ and
$L_2(W)=\gamma\cdot\Vec{d}$.
Note that $L_1(W)$ and $L_2(W)$ are equal to:
\begin{equation}\notag
L_1(W)=\#(\text{leaves})-2\#(\text{bridges}),\qquad
L_2(W)=\sum_{T:\text{ VEV tree}}n_{\text{root}}.
\end{equation}
${\mathcal H}$ is called the combined amplitude.

Then the partition function is expressed 
with 
${\rm Comb}_{\gamma}^{\bullet}
(\vec{\mu},\vec{\nu},\vec{\lambda})$
and
${\mathcal H}$:
\begin{prop}%
\begin{equation}\notag
{\mathcal Z}^{\gamma}_{\Vec{d}}(q)
=\sum_{\begin{subarray}{c}
(\Vec{\mu},\Vec{\nu},\Vec{\lambda});\\
\text{$r$-set of}\\
\text{degree $\Vec{d}$}
\end{subarray}}
\frac{1}
{z_{\Vec{\lambda}}z_{\Vec{\mu}}z_{\Vec{\nu}}}
\sum_{W\in {\rm Comb}^{\bullet}_{\gamma}
(\Vec{\mu},\Vec{\nu},\Vec{\lambda})}
{\mathcal H}(W).
\end{equation}
\end{prop}

\section{Free Energy}\label{feen}
In this section,
we obtain the free energy as 
the sum over
connected combined forests.

It is well-known,
for example in the calculation of the 
Feynman diagram,
that
if one takes the logarithm of
the sum of amplitudes over not necessarily connected  graphs,
then one obtains the sum over 
connected graphs.
More precisely, the statement is written as follows.
\begin{equation}\notag
\log \Bigg[
1+
\sum_{G\in \mathbb{G}^{\bullet}}
\frac{1}{|{\rm Aut}(G)|}\Psi(G)
\Bigg]
=
\sum_{G\in \mathbb{G}^{\circ}}
\frac{1}{|{\rm Aut}(G)|}\Psi(G)
\end{equation}
where $\mathbb{G}^{\bullet}$
is a set of graphs,
$\Psi$ is an amplitude on it
and $\mathbb{G}^{\circ}\subset\mathbb{G}^{\bullet}$
is the subset of connected graphs.
This formula is nothing but a variation of 
the exponential formula (\cite{stan}, chapter 5).
Therefore
we call 
it the {\em exponential formula}.

Of course, for the exponential formula to hold,
$\mathbb{G}^{\bullet}$ and $\Psi$ must satisfy
certain conditions.
The most important are 
that $\mathbb{G}^{\circ}$ generates
$\mathbb{G}^{\bullet}$ by the graph union operation
and that
$\Psi$ is multiplicative with respect to the graph union.
The rigorous formulation of the formula 
will appear in  appendix \ref{exf}.

In view of the exponential formula, 
we expect that 
the free energy is 
the sum over connected graphs.
In fact, this is true.
\begin{prop}\label{free}
%
\begin{equation}\notag
{\mathcal F}^{\gamma}_{\Vec{d}}(q)
=
\sum_{\begin{subarray}{c}
(\Vec{\mu},\Vec{\nu},\Vec{\lambda});\\
\text{$r$-set of}\\
\text{degree $\Vec{d}$}
\end{subarray}}
\frac{1}
{z_{\Vec{\lambda}}z_{\Vec{\mu}}z_{\Vec{\nu}}}
\sum_{W\in {\rm Comb}^{\circ}_{\gamma}
(\Vec{\mu},\Vec{\nu},\Vec{\lambda})}
{\mathcal H}(W).
\end{equation}
\end{prop}%
We delegate the proof to appendix \ref{pfe}.

\section{Proof of Main Theorem: Outline}\label{pt}
We sketch the 
proof of theorem \ref{mainthm} in this section.

The combined forest 
which is 
the same as a combined forest $W$ 
except that all the vertex-labels are multiplied by $k\in\mathbb{N}$,
is denote by
$W_{(k)}$.

We first rewrite ${\mathcal F}_{\Vec{d}}^{\gamma}(q)$ as follows.
\begin{equation}\notag
{\mathcal F}_{\Vec{d}}^{\gamma}(q)=
\sum_{k;k|\Vec{d}}
\sum_{\begin{subarray}{c}
(\Vec{\mu},\Vec{\nu},\Vec{\lambda});\\
\text{$r$-set of}\\
\text{degree $\Vec{d}/k$};\\
{\rm gcd}(\Vec{\mu},\Vec{\nu},\Vec{\lambda})=1
\end{subarray}}
\frac{1}
{z_{k\Vec{\lambda}}z_{k\Vec{\mu}}z_{k\Vec{\nu}}}
\sum_{W\in {\rm Comb}^{\circ}_{\gamma}
(\Vec{\mu},\Vec{\nu},\Vec{\lambda})}
{\mathcal H}(W_{(k)}).
\end{equation}
Then ${G}_{\Vec{d}}^{\gamma}(q)$ is equal to
\begin{equation}\notag
\begin{split}
{G}_{\Vec{d}}^{\gamma}(q)&=
\sum_{k;k|\Vec{d}}
\sum_{\begin{subarray}{c}
(\Vec{\mu},\Vec{\nu},\Vec{\lambda});\\
\text{$r$-set of}\\
\text{degree $\Vec{d}/k$};\\
{\rm gcd}(\Vec{\mu},\Vec{\nu},\Vec{\lambda})=1
\end{subarray}}
\sum_{k';k'|k}
\frac{k'}{k}
\mu
\Big(\frac{k}{k'}\Big)
\frac{1}
{z_{k'\Vec{\lambda}}z_{k'\Vec{\mu}}z_{k'\Vec{\nu}}}
\sum_{W\in {\rm Comb}^{\circ}_{\gamma}
(\Vec{\mu},\Vec{\nu},\Vec{\lambda})}
{\mathcal H}(W_{(k')})|_{q\to q^{k/k'}}
\\
&=\sum_{k;k|\Vec{d}}
\sum_{\begin{subarray}{c}
(\Vec{\mu},\Vec{\nu},\Vec{\lambda});\\
\text{$r$-set of}\\
\text{degree $\Vec{d}/k$};\\
{\rm gcd}(\Vec{\mu},\Vec{\nu},\Vec{\lambda})=1
\end{subarray}}
\frac{1}{k\#\aut(\Vec{\mu})\#\aut(\Vec{\nu})\#\aut(\Vec{\lambda})}
\sum_{W\in {\rm Comb}^{\circ}_{\gamma}(\Vec{\mu},\Vec{\nu},\Vec{\lambda})}
{\mathcal G}_k(W),
\end{split}
\end{equation}
where
\begin{equation}\notag
{\mathcal G}_k(W)=
\sum_{k':k'|k}
\mu\Big(\frac{k}{k'}\Big)
{{k'}^{-l(\Vec{\mu})-l(\Vec{\nu})-l(\Vec{\lambda})+1}}
{\mathcal H}(W_{(k')})\Big|_{q\to q^{k/k'}}.
\end{equation}
Therefore if we show
\begin{equation}\label{gt}
t\cdot{\mathcal G}_k(W)\in\mathbb{Q}[t]
\end{equation}
where $t=[1]^2$,
then the proof of the theorem will be finished.
This is the subject of the next section.

Since the entire proof is quite long,
let us briefly explain the main idea.
The amplitude ${\mathcal H}(W)$ turns out to 
be a polynomial in $t$ if $W$ has the cycle rank $>0$
and clearly $(\ref{gt})$ holds.
In the case the cycle rank is zero,
${\mathcal H}(W)$ has poles.
However, the poles of ${\mathcal H}(W)$ 
and ${\mathcal H}(W_{(k)})$
are related in some way and 
almost all of them cancel in ${\mathcal G}_k(W)$,
leaving at most a pole at $t=0$.
Let us describe the situation in the simplest case,
when %
$W$ does not contain a certain type of VEV trees
($V_{III}(W)=\emptyset$ in the notation of the next section),
and $k$ is odd.
Then
\begin{equation}\notag
{\mathcal H}(W_{(k)})=
\frac{g_W\cdot k^{(\Vec{\mu})+l(\Vec{\nu})+l(\Vec{\lambda})-1}}{t_k}
+\text{ polynomial in $t$ }.
\end{equation}
Here $t_k=[k]^2$ and
$g_W$ is a constant determined by $W$.
%
Therefore, if $k$ is odd,
\begin{equation}\notag
{\mathcal G}_k(W)
=
\frac{g_W}{t_k}\sum_{k';k'|k}\mu\Big(\frac{k}{k'}\Big)
+\text{polynomial in $t$}.
\end{equation}
%
(\ref{gt}) follows from the formula of M\"obius function
\begin{equation}\label{keyformula}
\sum_{k':k'|k}\mu\Big(\frac{k}{k'}\Big)
=\begin{cases}
1&(k=1)\\
0&(k>1,k\in\mathbb{N}).
\end{cases}
\end{equation}
So the factor $1/t_k$ in ${\mathcal G}_k(W)$ vanishes 
unless $k=1$. 
The proofs in other cases are more complicated, 
but similar.

\section{Pole Structure}\label{ps}
%
In this section, we first
study poles of the amplitude ${\mathcal B}(T)$
and 
then move to the study of 
${\mathcal H}(W)$.
In subsection \ref{PS1},
we see that ${\mathcal B}(T)$ is actually a 
function of $t=[1]^2$ and
study its poles.
Then we study
${\mathcal H}(W)$ 
in subsection \ref{PCA}.
The proof of the main theorem is completed
at the last of this section when we give a proof of
(\ref{gt}).
%

\subsection{${\mathcal B}(T)$}\label{PS1}
In this subsection,
we study poles of the amplitudes ${\mathcal B}(T)$
of a VEV tree $T$.
We first state the main proposition 
and check some examples.
Next we collect some lemmas necessary for the proof
and prove the proposition.
%
\subsubsection{Main Results}
We set 
\begin{equation}\notag
t=[1]^2,\qquad t_k=[k]^2 \quad(k\in\mathbb{N}).
\end{equation}
It is known that $t_k$ is a polynomial in $t$ with 
integer coefficients \cite{bp}
\begin{equation}\notag
t_k=\sum_{j=1}^k\frac{k}{j}
\begin{pmatrix}
j+k-1\\2j-1
\end{pmatrix}
t^j\qquad \in\mathbb{Z}[t].
\end{equation}

We define $m(T)$ of a VEV tree $T$ by
\begin{equation}\notag
m(T)={\rm gcd}(\mu,\nu)
\qquad
\text{ if } 
T\in \overline{{\rm Graph}}^{\circ}_a(\mu,\nu).
\end{equation}

It turns out that ${\mathcal B}(T)$ has poles at 
$t_{m(T)}=0$.
\begin{prop}\label{poleB}
Let $T$ be a VEV tree in 
$\overline{{\rm Graph}}_{a}^{\circ}(\mu,\nu)$.
\begin{enumerate}
\item %
There exists $g_{T}\in \mathbb{Z}$ and
$f_{T}(t)\in \mathbb{Z}[t]$ such that
\begin{equation}\notag
{\mathcal B}(T)=\begin{cases}
\frac{g_{T}}{t_{m(T)}}+f_{T}(t_{m(T)})& 
(\text{$m(T)$ odd or $n_{\text{root}}/m(T)$ even})
\\
\frac{g_{T}}{t_{m(T)}}
\Big(1+\frac{t_{m(T)/2}}{2}\Big)+f_{T}(t_{m(T)})& 
(\text{$m(T)$ even and $n_{\text{root}}/m(T)$ odd})
\end{cases}
\end{equation}
We remark that $n_{\text{root}}=ad$ where $d=|\mu|=|\nu|$.

\item %
Moreover,
\begin{equation}\notag
g_{T}=g_{T_{(0)}}\cdot {m(T)}^{l(\mu)+l(\nu)-1}.
\end{equation}
Here  $T_{(0)}$ is the VEV tree 
which is the 
same as $T$
except all the vertex-labels is multiplied by $1/m(T)$.
\end{enumerate}
\end{prop}%
We will give a proof for 
the case $m(T)=1$ first, then the case $m(T)>1$.

Before giving a proof,
let us check simple examples.
\newcommand{\dda}{
\unitlength .15cm
\begin{picture}(12,10)(0,0)
\thicklines
\put(1,7){$(d,0)$}
\put(7,7){$(-d,ad)$}
\put(5,5){\line(1,-1){2.5}}
\put(10,5){\line(-1,-1){2.5}}
\put(5,5){\circle*{1}}
\put(10,5){\circle*{1}}
\put(7.5,2.5){\circle*{1}}
\put(5,0){$(0,ad)$}
\end{picture}
}
\newcommand{\ddb}{
\unitlength .15cm
\begin{picture}(12,10)(0,0)
\thicklines
\put(1,7){$(1,0)$}
\put(7,7){$(-1,a)$}
\put(5,5){\line(1,-1){2.5}}
\put(10,5){\line(-1,-1){2.5}}
\put(5,5){\circle*{1}}
\put(10,5){\circle*{1}}
\put(7.5,2.5){\circle*{1}}
\put(5,0){$(0,a)$}
\end{picture}
}
Consider the VEV forest $T\in\overline{{\rm Graph}}^{\circ}_{a}((d),(d))$
($a\neq 0$)
below. $T$ has $m(T)=d$. 
\begin{equation}\notag
T=\raisebox{-.5cm}{\dda}
\qquad\qquad\qquad
T_{(0)}=\raisebox{-.5cm}{\ddb}
\end{equation}
If $a=1$ and $d=3$,
\begin{equation}\notag
\begin{split}
{\mathcal B}(T)
&=\frac{[9]}{[3]^3}
=\frac{1}{[3]^2}(q^{3}+1+q^{-3})
=\frac{t_3+3}{t_3},\qquad
g_T=3.
\\
{\mathcal B}(T_{(0)})
&=\frac{[1]}{[1]^3}=\frac{1}{t},\qquad
g_{T_{(0)}}=1.
\end{split}
\end{equation}
Here we have used the relation
\begin{equation}\notag
q^{k}+q^{-k}=t_k+2.
\end{equation}
If $a=3$ and $d=2$, 
\begin{equation}\notag
\begin{split}
{\mathcal B}(T)
&=\frac{[12]}{[6][2]^2}
=\frac{1}{[2]^2}(q^{3}+q^{-3})
=\frac{t_3+2}{t_2}
=\frac{t+2}{t_2}+t+2,
\quad
g_T=2.
\\
{\mathcal B}(T_{(0)})
&=\frac{[3]}{[3][1]^2}=\frac{1}{t},\qquad
g_{T_{(0)}}=1.
\end{split}
\end{equation}
Thus the proposition holds in these examples.

\subsubsection{Lemmas}\label{QN}
We collect some lemmas
on $q$-numbers necessary for the proof
of the proposition \ref{poleB}.
Here is the  list indicating where 
these lemmas are needed.

In the proof of the $k=1$ case,
lemma \ref{lemqnum}.4
plays the guiding role.
Lemmas \ref{G3} and \ref{G1} are used
to show that the third assumption of lemma \ref{lemqnum}.4.
In the proof of the $k>1$ case,
we use lemmas \ref{G3}.2 and \ref{G2}
together with the result of the $k=1$ case.
Lemmas \ref{iso} and \ref{lemqnum}.1,2,3 are
the preliminary  for the other lemmas.
We remark that the proof of lemma \ref{G2}
is very similar to the $k=1$ case.

Let us start from some notations.
For a (dummy) variable $x$,
the ring of polynomials with integer (rational) coefficients
is denoted by $\mathbb{Z}[x]$ ($\mathbb{Q}[x]$),
the ring of rational functions with rational coefficients
by $\mathbb{Q}(x)$.
The ring of Laurent polynomials with integer coefficients
is denoted by
$\mathbb{Z}[x,x^{-1}]$.
We define 
\begin{equation}\notag
\mathbb{Z}^+[x,x^{-1}]=\{f(x)\in \mathbb{Z}[x,x^{-1}]|
f(x)=f(x^{-1})\}.
\end{equation}
This is a subring of $\mathbb{Z}[x,x^{-1}]$.

We also define 
\begin{equation}\notag
{\mathcal L}[x]=\Big\{\frac{f_2(x)}{f_1(x)}\Big|
f_1(x),f_2(x)\in \mathbb{Z}[x],f_1(x):\text{ monic}
\Big\}.
\end{equation}
It is not difficult to see that this is a subring of
$\mathbb{Q}(x)$ (see \cite{peng}).

We set
\begin{equation}\notag
t=[1]^2,\qquad t_k=[k]^2,
\qquad
y=[1/2]^2.
\end{equation}
Note that $t=y(y+4)$.
%
\newcommand{\myapp}{<_y}
We define the relation $\myapp$ in ${\mathcal L}[y]$
by 
\begin{equation}\notag
f_1(y)\myapp f_2(y)\Leftrightarrow
{^\exists}f_3(y)\in\mathbb{Z}[y], f_3(0)=1,
\quad
f_1(y)=f_2(y)f_3(y).
\end{equation}
$f_1(y)\myapp f_2(y)$ implies that
$f_1(y)$ has poles at most where $f_2(y)$ does.
%
\begin{lemma}\label{iso}
\begin{equation}\notag
\mathbb{Z}^+[q,q^{-1}]
\cong \mathbb{Z}[t],
\qquad
\mathbb{Z}^+[q^{\frac{1}{2}},q^{-\frac{1}{2}}]
\cong \mathbb{Z}[y].
\end{equation}
\end{lemma}
\begin{proof}
We only show the
first isomorphism.
The proof of the 
second is the same if $q$ is replaced by 
$q^{\frac{1}{2}}$.

$\subseteq$:
An element $f(q)$ in $\mathbb{Z}^+[q,q^{-1}]$
is written as follows.
\begin{equation}\notag
f(q)=b_0+\sum_{i=1}^{m}
b_j(q^i+q^{-i})
=b_0+\sum_{i=1}^m b_j(t_{j}+2)
\qquad
(b_j\in \mathbb{Z}),
\end{equation}
where we have used the relation $ t_j= q^{j}+q^{-j}-2$.
Since $t_j\in\mathbb{Z}[t]$, $f(q)\in \mathbb{Z}[t]$.

$\supseteq$:
An element $f(t)$ in $\mathbb{Z}[t]$ is written as follows.
\begin{equation}\notag
f(t)=\sum_{i=0}^m
b_it^i
=\sum_{i=0}^m
b_i(q+q^{-1}-2)^i
\qquad
(b_j\in \mathbb{Z}).
\end{equation}
Thus $f(t)\in \mathbb{Z}^+[q,q^{-1}]$.
\end{proof}

\begin{lemma}\label{lemqnum}%
Let $a_1,\ldots,a_m,b_1,\ldots,b_n$ be integers.
\begin{enumerate}
\item 
$\prod_{i=1}^m[a_i]\in
\mathbb{Z}[y]\Leftrightarrow \text{ $m$ is  even}$.
%
%
\item
$\prod_{i=1}^m[a_i]\in
\mathbb{Z}[t]\Leftrightarrow 
\text{ both $m$ and $\sum_{i=1}^m a_i$ are even}$.
%
\item
$\prod_{i=1}^m[a_i]/\prod_{i=1}^n[b_i]\in
{\mathcal L}[y]\Leftrightarrow 
\text{ $m+n$ is  even}$.
\item
If
both $m+n$ and $ \sum_{i=1}^m a_i+\sum_{i=1}^n b_i $ are even
and
$\prod_{i=1}^m[a_i]/\prod_{i=1}^n[b_i]\in
\mathbb{Z}[q^{\frac{1}{2}},q^{-\frac{1}{2}}]$,
then 
$\prod_{i=1}^m[a_i]/\prod_{i=1}^n[b_i]\in
\mathbb{Z}[t].$
%
\end{enumerate}
\end{lemma}
\begin{proof}

1 and 2 follow immediately from the previous lemma.

3. If $m+n$ is odd,
$\prod[a_i]/\prod[b_j]$ is anti-symmetric under $q\to q^{-1}$.
Thus it can not be written with $y$.
If both $m$ and $n$ are even, then
$\prod[a_i],\prod[b_j]\in \mathbb{Z}[y]$.
If both $m$ and $n$ are odd,
then $[1]\prod[a_i],[1]\prod[b_j]\in \mathbb{Z}[y]$. 
Thus 3. is proved.

4.
The condition  that $m+n$ is even  means
that $\prod[a_i]/\prod[b_j]\in {\mathcal L}[y]$.
The assumption that it is a polynomial implies
that $l=\sum |a_i|-\sum |b_j|\geq 0$
and that it is written as follows.
\begin{equation}\notag
\frac{\prod[a_i]}{\prod[b_j]}
=\sum_{i=0}^l
c_i q^{\frac{l}{2}-i}
\qquad
(c_i\in\mathbb{Z},\quad c_i=c_{l-i}).
\end{equation}
The condition that $l$ is even implies that
$\prod[a_i]/\prod[b_j]\in
\mathbb{Z}^+[q,q^{-1}]$. 
Therefore 4 follows from the previous lemma.
\end{proof}

Recall that 
$[\lambda]=\prod_{i=1}^{l(\lambda)}[\lambda_i]$
where $\lambda$ is 
a partition.
\begin{lemma}\label{G3}
Let $k\in\mathbb{N}$ be a positive integer.
\begin{enumerate}
\item Let $a\in\mathbb{N}$ be a positive integer.
Then
$[ka]/[a]\in\mathbb{Z}[y]$. 
Moreover if $k$ is odd or $a$ is even, 
$[ka]/[a]\in\mathbb{Z}[t]$
and the constant term is $k$. 
If $k$ is even and $a$ is odd,
the remainder of $[ka]/[a]$ by $y(y+4)$
is $k(1+y/2)$.
\item
Let $\lambda$ be a partition. Then
$[k\lambda]/[\lambda]$ is written as follows.
\begin{equation}\notag
\frac{[k\lambda]}{[\lambda]}
=
\begin{cases}
k^{l(\lambda)}+t\times u_{k,\lambda}(t)
& (\text{$k$ odd or $|\lambda|$ even})
\\
k^{l(\lambda)}\Big(1+\frac{y}{2}\Big)
+t\times u_{k,\lambda}(y)
&(\text{$k$ even and $|\lambda|$ odd}).
\end{cases}
\end{equation}
Here, if $k$ is odd or $|\lambda|$ is even, 
$u_{k,\lambda}(t)\in \mathbb{Z}[t]$ is a polynomial of 
degree $(k-1)|\lambda|/2-1$,
if $k$ is even and $|\lambda|$ is odd,
$u_{k,\lambda}(y)\in\mathbb{Z}[y]$ is a polynomial of 
degree $(k-1)|\lambda|-2$.
\end{enumerate}
%
\end{lemma} 
%
Weaker statements of the lemma are
written with $\myapp$ in the following simple form:
\begin{equation}\notag
\frac{[ka]}{[a]}\myapp k,\qquad
\frac{[k\lambda]}{[\lambda]}\myapp k^{l(\lambda)}.
\end{equation}

\begin{proof}
1. Recall the formula $(x^k-1)/(x-1)=1+x+\cdots+x^{k-1}$.
Then it is easily calculated that 
\begin{equation}\notag
\frac{[ka]}{[a]}=
\sum_{i=0}^{k-1} q^{\frac{(k-1)a}{2}-i}.
\end{equation}
If $(k-1)a$ is even, $[ka]/[a]\in \mathbb{Z}^+[q,q^{-1}]$.
Thus it is in $\mathbb{Z}[t]$.
Since $t=0$ is equivalent to $q=1$,
the value at $t=0$ is the value at  $q=1$, which is $k$.

If $(k-1)a$ is odd, $[ka]/[a]\in \mathbb{Z}[y]$.
To compute the remainder by $y(y+4)$,
let us write it as $b_1+b_2y$.
We compare the values at $y=0$ and $y=-4$,
or equivalently at $q^{\frac{1}{2}}=1$ and 
$q^{\frac{1}{4}}=\sqrt{-1}$:
\begin{equation}\notag
b_1=k,\qquad
b_1-4b_2=-k.
\end{equation}
Therefore the first part is proved.

2. By the result of 1, 
$[k\lambda]/[\lambda]$ is in $\mathbb{Z}[y]$
with constant term $k^{l(\lambda)}$.

If $k$ is odd or $|\lambda|$ is even,
all the three conditions of lemma \ref{lemqnum}.4
are satisfied and $[k\lambda]/[\lambda]\in\mathbb{Z}[t]$.

If $k$ is even and $|\lambda|$ is odd,
there are odd number of odd parts in $\lambda$.
Take one odd part, say $\lambda_1$, and
denote by $\lambda'$ by the partition 
$\lambda\setminus \{\lambda_1\}$.
Then $[k\lambda']/[\lambda']\in\mathbb{Z}[t]$ with 
constant term $k^{l(\lambda)-1}$.
On the other hand,
$[k\lambda_1]/[\lambda_1]\in \mathbb{Z}[y]$ and the remainder
by $t=y(y+4)$ is $k(1+y/2)$.
Therefore the second part of the lemma follows.
\end{proof}%

\begin{lemma}\label{G1}
Let  $a,b,c\in{\mathbb{N}}$ be positive integers 
satisfying ${\rm gcd}(a,b,c)=1$. Then
\begin{equation}\notag
\frac{[{\rm lcm}(a,b,c)][{\rm gcd}(a,b)]
[{\rm gcd}(b,c)][{\rm gcd}(c,a)]
}{[a][b][c][1]}
\in \mathbb{Z}[t],
\end{equation}
and the constant term is $1$.
\end{lemma}%
With $\myapp$, (a weaker statement of) 
the lemma
are 
written as follows:
\begin{equation}\notag
\frac{[{\rm lcm}(a,b,c)]}{[a][b][c]}
\myapp \frac{[1]}{[{\rm gcd}(a,b)][{\rm gcd}(b,c)][{\rm gcd}(c,a)]}
.
\end{equation}
\begin{proof}
We check the three conditions to apply lemma \ref{lemqnum}.4.
The first condition is clearly satisfied.

We write $a,b,c$ as follows:
\begin{equation}\notag
a=pra',\qquad b=psb',\qquad c=qsc',
\end{equation}
where $p={\rm gcd}(a,b)$,
$s={\rm gcd}(b,c)$, $r={\rm gcd}(a,c)$.
Since ${\rm gcd}(a,b,c)=1$, any pair of 
$p,s,r,a',b',c'$ are prime to each other.

Then
\begin{equation}\notag
\begin{split}
l&:=psra'b'c'+p+s+r-pra'-psb'-src'-1\\
&=p(ra'-1)(sb'-1)+s(rc'-1)(sb'-1)-(p-1)(sb'-1)\equiv 0\mod2.
\end{split}
\end{equation}
Thus the second condition is also satisfied.

To check the third condition, we write:
\begin{equation}\notag
\begin{split}
&\frac{[{\rm lcm}(a,b,c)][{\rm gcd}(a,b)]
[{\rm gcd}(b,c)][{\rm gcd}(c,a)]
}{[a][b][c][1]}
=
\frac{[psra'b'c'][p][r][s]}{[pra'][psb'][src'][1]}
\\
=&q^{-l}
\underbrace{
\frac{(q^{psra'b'c'}-1)(q^p-1)(q^s-1)(q^r-1)}
{(q^{pra'}-1)(q^{psb'}-1)(q^{src'}-1)(q-1)}
}_{\star}.
\end{split}
\end{equation}
The numerator of $(\star)$
has a  zero of 4-th order at $q=1$,
double zeros at $q=e^{2\pi\sqrt{-1}j}$
$(j\in I_1)$
and
simple zeros at $q=e^{2\pi \sqrt{-1}j}$
$(j\in I_2)$
where
\begin{equation}\notag
\begin{split}
I_1&=\Big\{\frac{1}{p},\ldots,\frac{p-1}{p},
\frac{1}{s},\ldots,\frac{s-1}{s},
\frac{1}{r},\ldots,\frac{r-1}{r}
\Big\}
,
\\
I_2&=\Big\{\frac{1}{psra'b'c'},\ldots,
\frac{psra'b'c'-1}{psra'b'c}\Big\}
\setminus I_1.
\end{split}
\end{equation}
On the other hand,
the denominator has a  zero of 4-th order at $q=1$,
double zeros at $I_1$
and
simple zeros at $q=e^{2\pi\sqrt{-1}j}$
$(j\in I_3)$ 
where
\begin{equation}\notag
I_3
=\Big\{\frac{1}{pra'},\ldots,\frac{pra'-1}{pra'},
\frac{1}{psb'},\ldots,\frac{psb'-1}{psb'},
\frac{1}{src'},\ldots,\frac{src'-1}{src'}
\Big\}
\setminus I_1.
\end{equation}
Since $I_1\subseteq I_3$,
$(\star)$ is a polynomial in $q$.
Moreover, since
both the numerator and the denominator are elements of 
$\mathbb{Z}[q]$ and monic,
it is a polynomial with integer coefficients.
Therefore  the third condition is satisfied.
By lemma \ref{lemqnum}.4,
$[psra'b'c'][p][s][r]/[pra'][psb'][src'][1]\in\mathbb{Z}[t]$.

The constant term is equal to, by lemma \ref{G3}.1,
\begin{equation}\notag
\frac{[psra'b'c']}{[pra']}\Big|_{q=1}
\cdot 
\frac{[p]}{[psb']}\Big|_{q=1}
\cdot \frac{[s]}{[src']}\Big|_{q=1}
\cdot \frac{[r]}{[1]}\Big|_{q=1}
=
sb'c'\cdot \frac{1}{sb'}\cdot \frac{1}{rc'}\cdot r
=1.
\end{equation}
\end{proof}

\begin{lemma}\label{G2} 
Let $\lambda$ be a partition and $k\in \mathbb{N}$.
Assume they 
satisfy the following conditions.
1. $l(\lambda)\geq 2$.
2. ${\rm gcd}(k,\lambda\setminus \lambda_i)=1 
(1\leq {}^\forall i\leq l(\lambda))$.
3. $k$ is odd or $|\lambda|$ is even.
Then there exists $w_{k,\lambda}(t)\in \mathbb{Z}[t]$ and
\begin{equation}\notag
\frac{[k\lambda]}{[k]^2[\lambda]}=
\frac{k^{l(\lambda)-2}}{t}+w_{k,\lambda}(t).
\end{equation}
\end{lemma}

%

\begin{proof}
Define $p_i={\rm gcd}(k,\lambda_1,\ldots,\lambda_i)$ 
$(1\leq i\leq l(\lambda))$, $p_0=k$,
$k_i={\rm gcd}(k,\lambda_{i+1},\ldots,\lambda_{l(\lambda)})$
$(0\leq i\leq l(\lambda))$, $k_{l(\lambda)}=k$.
Note that ${\rm gcd}(\lambda_i,p_{i-1})=p_i$,
${\rm gcd}(\lambda_i,k_{i})=k_{i-1}$ and
${\rm gcd}(p_{i-1},k_i)=$
$p_{l(\lambda)}=k_0=1$ because of the second condition.

We write the LHS as follows.
\begin{equation}\notag
\frac{[k\lambda]}{[k]^2[\lambda]}
=\frac{1}{[k]^2}
\prod_{i=1}^{l(\lambda)}
\frac{[k\lambda_i]}{[{\rm lcm}(\lambda_i,p_{i-1},k_i)]}
\prod_{i=1}^{l(\lambda)}
\frac{[{\rm lcm} (\lambda_i,p_{i-1},k_i) ]}
{[\lambda_i][p_{i-1}][k_i]}
\prod_{i=1}^{l(\lambda)}
[p_{i-1}][k_i].
\end{equation}
By lemmas \ref{G1} and \ref{G3}.1,
\begin{equation}\notag
\begin{split}
(RHS)&
\myapp
\frac{1}{[k]^2}
\prod_{i=1}^{l(\lambda)}
\frac{1}{[p_i][k_{i-1}]}
\prod_{i=1}^{l(\lambda)}
[p_{i-1}][k_i]
\cdot
\prod_{i=1}^{l(\lambda)}
\frac{k\lambda_i}{{\rm lcm}(\lambda_i,p_{i-1},k_i)}
\\&
=\frac{1}{[1]^2}
\prod_{i=1}^{l(\lambda)}
\frac{k\lambda_i}{{\rm lcm}(\lambda_i,p_{i-1},k_i)}
=\frac{1}{[1]^2}\times(\text{positive integer})
.
\end{split}
\end{equation}

Therefore 
$[k\lambda][1]^2/[k]^2[\lambda]$ is in $\mathbb{Z}[y]$.
Applying lemma \ref{lemqnum}.4,
we find that 
$[k\lambda][1]^2/[k]^2[\lambda]\in \mathbb{Z}[t]$.
Its constant term is 
$[k\lambda]/[\lambda]|_{y=0}\cdot
[1]^2/[k]^2|_{y=0}=k^{l(\lambda)-2}
$.
\end{proof}

%
\subsubsection{Proof of Proposition \ref{poleB}: Case $m(T)=1$}%
We first show the case $m(T)=1$.

The statement actually
holds for a more general
VEV tree.
Let $\Vec{c}=(c_1,\ldots,c_l)$,
$\Vec{n}=(n_1,\ldots,n_l)$ 
be the pair satisfying the assumptions of
the subsection \ref{GD}
and ${\rm gcd}(\Vec{c},\Vec{n})=1$.
We extend
the definition of ${\mathcal B}(T)$ 
to $T\in \overline{{\rm Graph}}^{\circ}(\Vec{c},\Vec{n})$
as follows:
\begin{equation}\notag
{\mathcal B}(T)=\frac{{\mathcal A}(T)}{\prod_{i=1}^l[m_i]}
\end{equation}
where
$m_i={\rm gcd}(c_i,n_i)$ $(1\leq i\leq l)$.
We will show $t \cdot{\mathcal B}(T) \in\mathbb{Z}[t]$.
We give a proof of the case $|\Vec{n}|\neq 0$
and omit
the case $|\Vec{n}|=0$
since the proof is almost the same.

Our strategy of the proof is lemma \ref{lemqnum}.4.
Since
$[1]^2{\mathcal B}(T)$ is written by
the products of $q$-numbers:
\begin{equation}\notag 
[1]^2{\mathcal B}(T)=
\frac{[1]^2\prod_{v\in V_2(T)}[\zeta_v]}
{[|\Vec{n}|]\prod_{i=1}^l[m_i]},
\end{equation}
it is sufficient to show the followings:
%
\begin{equation}\notag
\begin{split}
1.\,&\#V_2(T)-l-1\equiv 0 \mod2;
\\
2.\,&\sum_{v\in V_2(T)}\zeta_v-|\Vec{n}|-\sum_{i=1}^lm_i\equiv0
\mod2;
\\
3.\,& [1]^2{\mathcal B}(T)\in \mathbb{Z}[y]
\end{split}
\end{equation}

1. The first condition is immediately checked since
$\#V_2(T)=l-1$.

2. 
Since $c_v=c_{L(v)}+c_{R(v)}$ and
$n_v=n_{L(v)}+n_{R(v)}$, 
\begin{equation}\notag
\zeta_v
=c_{L(v)}n_{L(v)}+c_{R(v)}n_{R(v)}-c_v n_v+2c_{L(v)}n_{R(v)}.
\end{equation}
Then
\begin{equation}\notag
\begin{split}
&\sum_{v\in V_2(F)}\zeta_v-|\Vec{n}|-\sum_{i=1}^lm_i
\mod2
\\
\equiv
&\sum_{i=1}^l c_in_i-
|\Vec{n}|-\sum_{i=1}^lm_i-|\Vec{c}|
\mod2
\qquad(\because|\Vec{c}|=0)
\\
\equiv&
\sum_{i=1}^l\big[m_i^2(c_i'-1)(n_i'-1)+
m_i(m_i-1)(c_i'+n_i')-m_i(m_i-1)\big]\mod2
\\
\equiv& 0\mod 2
\end{split}
\end{equation}
where $c_i',n_i'$ are defined by 
$c_i=m_ic_i'$ and $n_i=m_in_i'$ for $1\leq i\leq l$.

3.
Let us  introduce some notations.
For a vertex $v$,
\begin{equation}\notag
\begin{split}
\Vec{c}_v&=
\{c_i\in \Vec{c}| \text{ the $i$-th leaf is a descendent of $v$}\},\\
\Vec{n}_v&=\{n_i\in \Vec{n}| \text{ the $i$-th leaf is a 
descendent of $v$}\},
\\
p_v&={\rm gcd}(\Vec{c}_v,\Vec{n_v}),\qquad
k_v={\rm gcd}(\Vec{c}\setminus \Vec{c}_v,
\Vec{n}\setminus\Vec{n}_v,|\Vec{n}|).
\end{split}
\end{equation}
$p_v,k_v$'s satisfy the following relations. 
(For $a,b\in\mathbb{N}$, $a|b$ means $b$ is divisible by $a$.)
\begin{equation}\notag
\begin{split}
&p_v,k_v|\,{\rm gcd}(c_v,n_v),\qquad
p_{L(v)},p_{R(v)},k_v|\zeta_v,
\\
&{\rm gcd}(p_{R(v)},k_v)=k_{L(v)},\qquad
{\rm gcd}(p_{L(v)},k_v)=k_{R(v)},
\\
&p_{\text{$i$-th leaf}}=m_i,\qquad
k_{\text{root}}=|\Vec{n}|,\qquad
{\rm gcd}(p_v,k_v)=k_{\text{leaf}}=p_{\text{root}}=1.
\end{split}
\end{equation}
Note that the last identities are the consequence of 
the assumption
${\rm gcd}(\Vec{c},\Vec{v})=1$.

Now we decompose the amplitude ${\mathcal B}(T)$ 
as follows.
\begin{equation}\notag
\begin{split}
{\mathcal B}(T)
&=
\underbrace{
\prod_{v\in V_2(T)}
\frac{[\zeta_v]}
{[{\rm lcm}(k_v,p_{L(v)},p_{R(v)})]}
}_{(a)}
\underbrace{
\prod_{v\in V_2(T)}
\frac{[{\rm lcm}(k_v,p_{L(v)},p_{R(v)})]}
{[k_v][p_{L(v)}][p_{R(v)}]}
}_{(b)}
\\&\times
\underbrace{
\prod_{v\in V_2(T)}
{[k_v][p_{L(v)}][p_{R(v)}]}
}_{(c)}
\cdot
\underbrace{
\frac{1}{[|\Vec{n}|]\prod_{i=1}^l[m_i]}
}_{(d)}.
\end{split}
\end{equation}
Since $k_v,p_{L(v)},p_{R(v)} $ divide $\zeta_v$,
so does ${\rm lcm}( k_v,p_{L(v)},p_{R(v)})$.
By  lemma \ref{G3}.1,
\begin{equation}\notag
(a)\myapp 
\prod_{v\in V_2(T)}\frac{\zeta_v}{
 {\rm lcm}( k_v,p_{L(v)},p_{R(v)})}=:g_{T}\in\mathbb{Z}.
\end{equation}

We apply lemma \ref{G1} to $(b)$.
Note that the assumption is satisfied
since ${\rm gcd}(k_v,p_{L(v)},p_{R(v)})=1$.
\begin{equation}\notag
\begin{split}
(b)&\myapp \prod_{v\in V_2(T)}
\frac{[1]}{[k_{L(v)}][k_{R(v)}][p_v]}
\\
&=
[1]^{\#V_2(T)}
\cdot \frac{1}{[p_{\text{root}}]}
\cdot
\prod_{v:\text{leaf}} \frac{1}{[k_v]} 
\cdot
\prod_{\begin{subarray}{c}v\in V_2(T)\\
            v\neq \text{root}\end{subarray}}
\frac{1}{[p_v][k_v]}
\\
&=
\frac{1}{[1]^{2}}
\prod_{\begin{subarray}{l}v\in V_2(T),\\v\neq \text{root}
\end{subarray}}
\frac{1}{[p_v][k_v]}.
\end{split}
\end{equation}
We have used $[k_{\text{leaf}}]=[p_{\text{root}}]=[1]$
and $\#V_2(T)=\#(\text{leaves})-1$.

The factor $(c)$ is rewritten as follows:
\begin{equation}\notag
\begin{split}
(c)&=
[k_{\text{root}}]\cdot
\prod_{v:\text{leaf}}[p_v] 
\cdot
\prod_{\begin{subarray}{c}v\in V_2(T)\\
 v\neq \text{root}\end{subarray}}
{[p_v][k_v]}
=
[|\Vec{n}|]\cdot
\prod_{i=1}^l[m_i]
\prod_{\begin{subarray}{c}v\in V_2(T),\\
v\neq \text{root}
\end{subarray}}
[p_v][k_v].
\end{split}
\end{equation}

Thus all the factors except $1/[1]^2$ cancel  and we obtain
\begin{equation}\notag
{\mathcal B}(T)\myapp \frac{g_T}{[1]^2}
.
\end{equation}
Thus we find that $[1]^2 {\mathcal B}(T)\in \mathbb{Z}[y]$.
All the assumptions of lemma \ref{lemqnum}.4 are checked
and
the $m(T)=1$ case is proved.

\subsubsection{Proof of Proposition \ref{poleB}: Case $m(T)>1$}
Next we give a proof in the case $m(T)>1$.
Here $T$ is a VEV tree in 
$\overline{{\rm Graph}}_{a}^{\circ}(\mu,\nu)$
with $|\mu|=|\nu|=d$.
For simplicity of the writing, let us set 
$k=m(T)={\rm gcd}(\mu,\nu)$.

${\mathcal B}(T)$ is written in
terms of ${\mathcal B}(T_{(0)})$ as follows.
\begin{equation}\notag
{\mathcal B}(T)=
\underbrace{
\frac{[k\mu'][k\nu']}
{[k]^2[\mu'][\nu']}\Big|_{q\to q^k}
}_{(a')}
\cdot
\underbrace{
\frac{[kad']}{[ad']}\Big|_{q\to q^k}
}_{(b')}
\times 
\underbrace{
\Big([1]^2{\mathcal B}(T_{(0)})\big)\Big|_{q\to q^{k^2}}}_{(c')}
\end{equation}
where $\mu',\nu'$ and 
$d'$ are defined by $\mu=k\mu',\nu=k\nu'$ and 
$d=kd'$.
%
%

We will apply lemma \ref{G2} to the factor
$(a')$ by 
substituting $\lambda'\cup\mu'$ into $\lambda$.
Let us check the assumptions.
The first and the third are clearly satisfied 
since
$l(\mu')+l(\nu')\geq 2$ and $|\mu'|+|\nu'|=2d'$.
The second assumption is also satisfied since
if we assume otherwise,
then it contradicts with $|\mu'|=|\nu'|$.
Therefore by the lemma,
there exists $w_{k, \mu'\cup\nu'}(t)\in\mathbb{Z}[t]$
such that
\begin{equation}\notag
(a')=\frac{k^{l(\mu)+l(\nu)-2}}{t_k}+
w_{k,\mu'\cup\nu'}(t_k).
\end{equation}

By lemma \ref{G3},
the factor $(b')$ is written as follows.
\begin{equation}\notag
(b')=\frac{[kad']}{[ad']}\Big|_{q\to q^k}
=\begin{cases}
k+t_k \times u_{k,(d')}(t_k)
& (\text{$k$ odd or $ad'$ even})
\\
k\Big(1+\frac{t_{{k}/{2}}}{2}\Big)
+t_k\times u_{k,(d')} (t_{{k}/{2}})
&(\text{$k$ even and $ad'$ odd})
\end{cases}
\end{equation}
where, in both cases, 
$u_{k,(d')}$ is a polynomial in one variable with
integer  coefficients.

From the result of the $m(T)=1$ case, 
$(c')\in\mathbb{Z}[t_{k^2}]\subset\mathbb{Z}[t_k]$
with the constant term $g_{T_{(0)}}$.
Thus the $m(T)>1$ case is proved.
This completes the proof of the proposition.

\subsection{Poles of Combined Amplitudes}\label{PCA}
We study the pole of 
the combined amplitude and 
finish the proof of the main theorem
by showing (\ref{gt}): 
$t\,{\mathcal G}_k(W)\in\mathbb{Q}[t]$.

As it turns out,
the pole structure of ${\mathcal H}(W)$
of $W\in {\rm Comb}_{\gamma}^{\circ}
(\Vec{\mu},\Vec{\nu},\Vec{\lambda})$
depends on the cycle rank $\beta(W)$
and $k={\rm gcd}(\Vec{\mu},\Vec{\nu},\Vec{\lambda})$.
We first see that
if $W$ has the cycle rank $>0$,
then ${\mathcal H}(W)$ is a polynomial in $t$
and if $W$ has the cycle rank zero,
${\mathcal H}(W)$ has a simple pole at $t_k=0$.
Then we prove the more detailed pole structure in  the latter case.
We also prove (\ref{gt}).
%

\subsubsection{}

\begin{prop}\label{combpole1}
Let $(\Vec{\mu},\Vec{\nu},\Vec{\lambda})$ be
an $r$-set
and $k={\rm gcd}(\Vec{\mu},\Vec{\nu},\Vec{\lambda})$.
Let  $W\in {\rm Comb}_{\gamma}^{\circ}
(\Vec{\mu},\Vec{\nu},\Vec{\lambda})$
be a connected combined forest.
\begin{enumerate}
\item
If $\beta(W)=0$, then
$t_k {\mathcal H}(W)\in\mathbb{Z}[t]$.
\item %
%
If $\beta(W)>1$, then ${\mathcal H}(W)\in\mathbb{Z}[t]$.
\end{enumerate}
\end{prop}
\begin{proof}
We write the combined amplitude as follows.
\begin{equation}\notag
{\mathcal H}(W)=(-1)^{L_1(W)+L_2(W)}
\prod_{T:\text{ VEV tree in $W$}}
{\mathcal B}(T)\cdot \prod_{b:\text{ bridges in $W$}}
t_{h(b)}.
\end{equation}
Since every ${\mathcal B}(T)$ has poles at $t_{m(T)}=0$,
the first factor has poles at
\begin{equation}\notag
\{t_{m(T)}=0|\text{ $T$: VEV tree in $W$}\}.
\end{equation}
On the other hand the second factor has zeros at
\begin{equation}\notag
\{t_{h(b)}=0|\text{ $b$: bridges in $W$}\}.
\end{equation} 
Do the poles cancel with the zeros ?
Note that if a bridge $b$ is incident to a VEV tree $T$,
then the zeros of $t_{h(b)}$ cancel the poles of
$t_{m(T)}$ because
$m(T)$ is a divisor of  ${h(b)}$.
Therefore the question is whether it is possible to couple every
VEV tree in $W$ with at least one of its incident bridges.

The answer is yes if $\beta(W)>0$.
Therefore ${\mathcal H}(W)$ is a polynomial in $t$.
As an example, consider the first graph in  figure \ref{excomb1}.
The cycle rank is two in this case. 
Let $T_1,T_2,T_3$ denote the VEV trees and 
$b_1,b_2,b_3$ the bridges as depicted below.
(For simplicity, we contract each VEV tree to a vertex.)
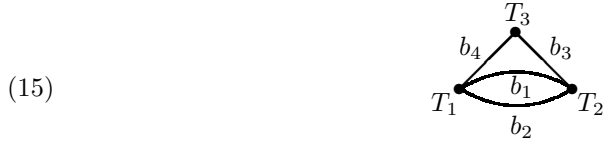
\newcommand{\fstgraph}{
\unitlength .15cm
\begin{picture}(15,10)(-5,-4)
\thicklines
\put(0,0){\line(1,1){5}}
\put(10,0){\line(-1,1){5}}
\qbezier(0,0)(5,3)(10,0)
\qbezier(0,0)(5,-3)(10,0)
\put(0,0){\circle*{1}}
\put(10,0){\circle*{1}}
\put(5,5){\circle*{1}}
\put(-2.5,-2){$T_1$}
\put(10.5,-2){$T_2$}
\put(4,6){$T_3$}
\put(4.5,-.5){$b_1$}
\put(4.5,-4){$b_2$}
\put(8,3){$b_3$}
\put(0,3){$b_4$}
\end{picture}
}
\begin{equation}\label{fst}
\raisebox{-.5cm}{
\fstgraph}
\end{equation}
If we take the coupling
$(T_3,b_3)$, $(T_2,b_1,b_2)$, $(T_1,b_4)$,
we could easily see that 
the pole of ${\mathcal B}(T_3)$ is canceled by 
the zeros of $t_{h(b_3)}$ etc.

\begin{figure}[th]
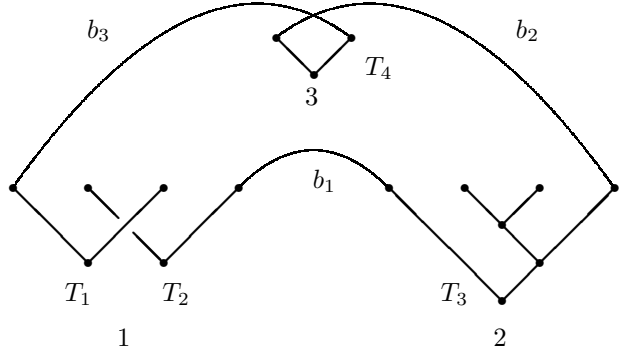

\begin{equation}\notag
\!\!\!\!\!\!\!\!\!\!
\!\!\!\!\!\!\!\!\!\!
\!\!\!\!\!\!\!\!\!\!
\!\!\!\!\!\!\!\!\!\!
\!\!\!\!\!\!\!\!\!\!
\!\!\!\!\!\!\!\!\!\!
\!\!\!\!\!\!\!\!\!\!
\!\!\!\!\!\!\!\!\!\!
\ancomb{\treeb}{\treea}
\end{equation}
\caption{A combined forest in ${\rm Comb}_{\gamma}^{\bullet}
(\vec{\mu},\vec{\nu},\vec{\lambda})$
where
$r=3$, $\gamma_i+2\neq 0$ and
$\lambda^1=(1)=\lambda^3=\lambda^2=(1)$,
$\mu^1=(1)=\mu^2=\nu^1=\nu^2, \mu^3=\nu^3=\emptyset$.
The cycle rank is zero.
}\label{excomb2}
\end{figure}
If the cycle rank $\beta(W)$ is zero,
the answer to the question is no.
However, if we pick one VEV tree $T$,
then it is possible to couple every other VEV tree with 
one of its incident bridges.
As an example, consider the combined forest in the figure 
\ref{excomb2}.
For simplicity, 
we contract VEV trees to a vertex as before (see below).
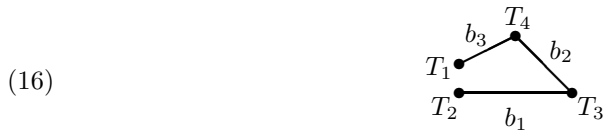
\newcommand{\scdgraph}{
\unitlength .15cm
\begin{picture}(15,10)(-5,-3)
\thicklines
\put(0,2.5){\line(2,1){5}}
\put(10,0){\line(-1,1){5}}
\put(0,0){\line(1,0){10}}
\put(0,0){\circle*{1}}
\put(0,2.5){\circle*{1}}
\put(10,0){\circle*{1}}
\put(5,5){\circle*{1}}
\put(-2.5,-2){$T_2$}
\put(10.5,-2){$T_3$}
\put(4,6){$T_4$}
\put(-3,1.5){$T_1$}
\put(4,-3){$b_1$}
\put(8,3){$b_2$}
\put(0.5,4.5){$b_3$}
\end{picture}
}
\begin{equation}\label{scd}
\raisebox{-.5cm}{
\scdgraph
}
\end{equation}
If we pick $T_4$ 
and choose the pairing
%
$(T_1,b_3)$, $(T_2,b_1)$, $(T_3,b_2)$,
then we could easily see that
the pole of ${\mathcal B}(T_1)$ cancel with the zeros
of $t_{h(b_3)}$, etc,
and that
only the pole of ${\mathcal B}(T_4)$ remains.
Moreover, the choice of a VEV tree $T$ is arbitrary.
Therefore ${\mathcal H}(W)$ has poles at most
\begin{equation}\notag
\begin{split}
\bigcap_{T:\text{ VEV tree in $W$}}
\{t_{m(T)}=0\}
&=
\{t_{{\rm gcd}(m(T),{T :\text{VEV tree in $W$}}))}=0\}
\\
&=\{t_{({\rm gcd}(\Vec{\mu},\Vec{\nu},\Vec{\lambda}))}=0\}=
\{t_k=0\}.
\end{split}
\end{equation}
Therefore
$t_k {\mathcal H}(W)$ is a polynomial in $t$.

The integrality of the coefficients follows from 
lemma \ref{lemqnum}.4.
\end{proof}

\subsubsection{}
We study the property of the case $\beta(W)=0$ in more details.
%
\newcommand{\mysimz}{\sim_{\mathbb{Z}}}%
\newcommand{\mysimq}{\sim_{\mathbb{Q}}}%
For convenience, we define the following two relations. 
\begin{itemize}
\item
The relation
$\mysimz$ in the ring ${\mathcal L}[t]$ is defined
by
\begin{equation}\notag
u_1(t)\mysimz u_2(t)\Leftrightarrow
u_1(t)-u_2(t)\in\mathbb{Z}[t].
\end{equation}
\item
The relation
$\mysimq$ in the ring $\mathbb{Q}(t)$ is defined
by
\begin{equation}\notag
u_1(t)\mysimq u_2(t)\Leftrightarrow
u_1(t)-u_2(t)\in\mathbb{Q}[t].
\end{equation}
\end{itemize}
Both $u_1(t)\mysimz u_2(t)$ and $u_1(t)\mysimq u_2(t)$ imply that
$u_1(t)$ and $u_2(t)$ has the same poles.

\begin{prop}
\label{lastprop}
Let $(\Vec{\mu},\Vec{\nu},\Vec{\lambda})$ be
an $r$-set
such that 
${\rm gcd}(\Vec{\mu},\Vec{\nu},\Vec{\lambda})=1$.
Let $W\in {\rm Comb}_{\gamma}^{\circ}
(\Vec{\mu},\Vec{\nu},\Vec{\lambda})$
be a connected combined forest with 
the cycle rank $\beta(W)=0$.
\begin{enumerate}
\item 
If $k$ is an odd positive integer,
\begin{equation}\notag
{\mathcal H}(W_{(k)})\mysimz
k^{l(\Vec{\mu})+l(\Vec{\nu})+l(\Vec{\lambda})-1}
{\mathcal H}(W)|_{t\to t_k}.
\end{equation}
\item If $k$ is an even positive integer,
\begin{equation}\notag
{\mathcal H}(W_{(k)})\mysimz
(-1)^{L_2(W)}
k^{l(\Vec{\mu})+l(\Vec{\nu})+l(\Vec{\lambda})-1}
{\mathcal H}(W)|_{t\to t_k}
\prod_{T\in VT_I(W)}
\Big(1+\frac{t_{m(T)k/2}}{2}\Big).
\end{equation}
$VT_I(W)$ denotes the set of VEV trees in $W$
such that both $m(T)$ and $n_{\text{root}}$ are odd.
%
\item %
\begin{equation}\notag
{\mathcal G}_k(W)\mysimq
\begin{cases}
0& (k>2)\\
\text{const.}/t &(k=1,2).
\end{cases}
\end{equation}
\end{enumerate}
\end{prop}%
As we have seen in  section \ref{pt},
the main theorem \ref{mainthm} follows directly from
the third statement.

Before giving a proof,
let us check the proposition in simple examples.
Consider the combined forest
in  figure \ref{excomb2} (and also see (\ref{scd}))
with $\gamma_1=\gamma_2=\gamma_3=-1$.
In this case
$l(\Vec{\mu})+l(\Vec{\nu})+l(\Vec{\lambda})=7$.
VEV trees of type I are $T_1,T_2,T_4$.
All VEV tree have $m(T_i)=1$.
The amplitude of $W$ is equal to
\begin{equation}\notag
{\mathcal H}(W)=\frac{1}{[1]^2}=\frac{1}{t}.
\end{equation}
As a check for the odd $k$ case, take $k=3$.
The amplitude of $W_{(3)}$ is
\begin{equation}\notag
\frac{3^6}{t_3}+
3^7+3^5 11 t_3+
3^3 5 \cdot 13 {t_3}^2+
3^3 5^2 {t_3}^3
+
3^2 17 {t_3}^4+
19{t_3}^5+t_3.
\end{equation}
As an example for the even $k$, we take $k=2$.
The amplitude of $W_{(2)}$ is
\begin{equation}\notag
\frac{2^6(1+\frac{t}{2})^2}{t_2}
+32+96t+86t^2+41t^3+10t^4+t^5.
\end{equation}
Therefore the first and the second 
of the proposition hold.
It is easy to see that the third statement also holds 
in these cases.

\begin{proof}%
1 and 2.
Since the pole structure of
${\mathcal B}(T)$
depends on $m(T)$ and $n_{\text{root}}$ (proposition \ref{poleB}),
we separate VEV trees in $W$ into three types:
\begin{equation}\notag
\begin{split}
VT_I(W)&=
\{T:\text{VEV tree}| 
\text{$m(T)$ odd and $n_{\text{root}}/m(T)$ odd}\},
\\
VT_{II}(W)
&=\{T:\text{VEV tree}|\text{$n_{\text{root}}/m(T)$ even}\},
\\
VT_{III}(W)&=
\{T:\text{VEV tree}| \text{$m(T)$ even and 
$n_{\text{root}}/m(T)$ odd}\}.
\end{split}
\end{equation}
The combined amplitude ${\mathcal H}(W)$ is written as follows.
\begin{equation}\label{une}
{\mathcal H}(W)
\mysimz
(-1)^{L_1(W)+L_2(W)}
\prod_{b:\text{bridge}}t_{h(b)}\cdot
\prod_{T: \text{VEV tree}}
\frac{g_T}{t_{m(T)}}
\cdot
\prod_{T\in V_{III}(W)}
\Big(1+\frac{t_{m(T)/2}}{2}\Big).
\end{equation}

Consider the amplitude of $W_{(k)}$.
If $k$ is odd,
\begin{equation}\notag
VT_{I}(W_{(k)})\cong VT_{I}(W),\qquad
VT_{II}(W_{(k)})\cong VT_{II}(W),\qquad
VT_{III}(W_{(k)})\cong VT_{III}(W).
\end{equation}
If $k$ is even,
\begin{equation}\notag
VT_{I}(W_{(k)})=\emptyset,\quad
VT_{II}(W_{(k)})\cong VT_{II}(W),\quad
VT_{III}(W_{(k)})\cong VT_{I}(W)\cup VT_{III}(W).
\end{equation}
Thus we have
\begin{equation}\label{deux}
\begin{split}
{\mathcal H}(W_{(k)})
&\mysimz
(-1)^{L_1(W)+kL_2(W)}
k^{\#(\text{leaves})-\#(\text{trees})}
\\&\times
\prod_{b:\text{bridge in $W$}}t_{kh(b)}\cdot
\prod_{T: \text{VEV tree in $W$}}
\frac{g_T}{t_{km(T)}}
\\&\times
\begin{cases}
\prod_{T\in V_{III}(W)}
\Big(1+\frac{t_{km(T)/2}}{2}\Big)& (k\text{:odd})
\\
\prod_{T\in V_{I}(W)\cup V_{III}(W)}
\Big(1+\frac{t_{km(T)/2}}{2}\Big)& (k\text{:even})
\end{cases}
\end{split}
\end{equation}
Given that $\beta(W)=\#(\text{bridges})-\#(\text{trees})+1=0$,
$\#(\text{bridges})=l(\Vec{\lambda})$
and that
$\#(\text{leaves})=l(\Vec{\mu})+l(\Vec{\nu})+2l(\Vec{\lambda})$,
\begin{equation}\notag
k^{\#(\text{leaves})-\#(\text{trees})}=
k^{l(\Vec{\mu})+l(\Vec{\nu})+l(\Vec{\lambda})-1}.
\end{equation}
Comparing (\ref{une}) and (\ref{deux}), 
we obtain
the first and the second statements.

3. 
In the proof of the third statement,
the key formula 
is the formula for the M\"obius function
(\ref{keyformula}).

(i) $k$ is odd. Since the divisors of $k$ are all odd,
\begin{equation}\notag
{\mathcal G}_k(W)\mysimq
{\mathcal H}(W)|_{q\to q^k}
\sum_{k';k'|k}
\mu\Big(\frac{k}{k'}\Big)
=\begin{cases}
{\mathcal H}(W)& (k=1)\\
0 &(k>1)
\end{cases}
\end{equation}
%
Therefore this case follows from the previous proposition.

(ii) $k$ is divisible by $4$.
If we write $k=4k_0$ $(k_0\geq 1)$,
every divisor $k'$ of $k$ is written as 
$4k'',2k'',k''$ with a divisor $k''$ of $k_0$.
For $k'=k''$,
the M\"obius function
$\mu(k/k')$ is zero since $k/k'$ contains $2^2$.
For $k'=4k'',2k''$,
we could write the result of proposition \ref{lastprop}.2 as 
\begin{equation}\notag
{\mathcal H}(W_{(k')})\mysimq 
{\mathcal H}(W_{(2)})|_{q\to q^{k'/2}}
\Big(\frac{k'}{2}\Big)^{l(\Vec{\mu})+l(\Vec{\nu})+l(\Vec{\lambda})-1}
\end{equation}
Therefore
\begin{equation}\notag
{\mathcal G}_k(W)\mysimq
\big(2^{-l(\Vec{\mu})-l(\Vec{\nu})-l(\Vec{\lambda})+1 }
{\mathcal H}(W_{(2)})|_{q\to q^{k/2}}
\sum_{k';k'|k}
\mu\Big(\frac{k}{k'}\Big)
=0.
\end{equation}
This proves the case $4|k$.

(iii)
$k=2k_0$ where $k_0>1$ is an odd positive integer.
Note that any divisor of $k$ is 
written as either $k'$ or $2k'$ with
a divisor $k'$ of $k_0$.
Therefore
\begin{equation}\notag
{\mathcal G}_k(W)\mysimq
\big(2^{-l(\Vec{\mu})-l(\Vec{\nu})-l(\Vec{\lambda})+1 }
{\mathcal H}(W_{(2)})|_{q\to q^{k/2}}
-{\mathcal H}(W)|_{q\to q^k}
\big)
\sum_{k';k'|k_0}
\mu\Big(\frac{k}{k'}\Big)
=0.
\end{equation}

(iv) $k=2$.
\begin{equation}\notag
\begin{split}
{\mathcal G}_2(W)
&\mysimq
-{\mathcal H}(W)|_{t\to t(t+4)}
\Big(1+
(-1)^{L_2(W)}\prod_{T\in VT_{I}(W)}
(1+\frac{t_{m(T)}}{2})
\Big)
\end{split}
\end{equation}
Since $\#VT_{I}(W)\equiv
\sum_{T}n_{\text{root}}\equiv L_2(W)\mod 2$,
the $k=2$ case follows from the next lemma.
\begin{lemma}
Let $m_1,m_2,\ldots,m_{l}$ be 
positive odd integers.
Then 
\begin{equation}\notag
\frac{(-1)^{l+1}+\prod_{i=1}^{l}(1+\frac{t_{m_i}}{2})}{t+4}\in
\mathbb{Q}[t].
\end{equation}
\end{lemma}
\begin{proof}
To show the lemma, we only have to show that
the denominator has zero at $t=-4$.
At $t=-4$, or equivalently at $q^{\frac{1}{2}}=\sqrt{-1}$,
$t_{m_i}=(q^{\frac{1}{2}})^{2m_i}+
(q^{\frac{1}{2}})^{-2m_i}-2=-4$.
Thus the denominator is zero at $t=-4$ and the lemma is proved.
\end{proof}

This completes the proof of the proposition \ref{lastprop}.
\end{proof}%

\appendix%
\section{Exponential Formula}\label{exf}
In this appendix,
we formulate the exponential formula
used 
in section \ref{feen}.

First, we present the conditions  for a graph set.
\begin{definition} \label{Gset}
A {\em G-set} is  a set of (labeled) graphs
$\mathbb{G}^{\bullet}$ 
together with a family of subsets,
$\{\mathbb{G}_{d}^{\bullet}\}_{d\geq 1}$,
satisfying the following conditions.
\begin{enumerate}
\item
$\mathbb{G}^{\bullet}=\coprod_{d\geq 1}\mathbb{G}_d^{\bullet}$.
\item \label{nonull}
If $G\in\mathbb{G}_d^{\bullet}$ and 
$G'\in\mathbb{G}_{d'}^{\bullet}$, then
$G\cup G'\in\mathbb{G}_{d+d'}^{\bullet}$.
%
\item \label{generation}
Every connected component $G_i$ of 
$G\in\mathbb{G}^{\bullet}$ is also an element of
$\mathbb{G}^{\bullet}$.
\item \label{nonisomorphic}
If $G$ and $G'$ are two  elements in $\mathbb{G}^{\bullet}$,
then $G$ and $G'$ are not isomorphic.
\item \label{finiteness}
Every $\mathbb{G}_d^{\bullet}$ $(d\geq 1)$ is a finite set.
\end{enumerate}
If $G\in\mathbb{G}^{\bullet}_d$,
the {\em degree} of $G$ is $d$ and is denoted by $d(G)$. 
\end{definition}

\begin{remark}\label{gen}
The conditions \ref{nonull} and \ref{generation}
together with the condition \ref{nonisomorphic}
imply that
the decomposition of a graph into connected components:
\begin{equation}\notag
G=G_1\cup\ldots\cup G_l\mapsto
\{G_1,\ldots, G_l\}
\end{equation}
gives a one-to-one correspondence
\begin{equation}\notag
\mathbb{G}_{(l)}^{\bullet}\stackrel{\sim}{\longrightarrow}
SP^l\mathbb{G}^{\circ},
\end{equation}
where
$\mathbb{G}^{\bullet}_{(l)}\subset\mathbb{G}^{\bullet}$ 
is the subset 
consisting of graphs with $l$ connected components
and
$SP^l\mathbb{G}^{\circ}$
is
the $l$-th symmetric product 
of $\mathbb{G}^{\circ}$,
the subset consisting of all connected graphs.
\end{remark}

Next we set the condition for an amplitude.
We recall the definition of 
a graded ring. 
A graded ring is a commutative associative ring $R$ 
together with a family $(R_d)_{d\geq 0}$
of subgroups of the additive group of $R$, such that
$R=\oplus_{d\geq 0}R_d$ and $R_{d}\cdot R_{d'}\subseteq R_{d+d'}$
for all $d,d'\geq 0$.

\begin{definition}\label{GAtriple}
Let $(\mathbb{G}^{\bullet},$R$,\Psi)$
be a triple of a G-set, a 
graded ring 
and 
a grading preserving map $\Psi:\mathbb{G}\to R$.
It is a 
{\em triple of graph amplitude}
(GA triple)
if  the map $\Psi$ is multiplicative with respect to the graph union:
\begin{equation}\notag
\Psi(G\cup G')=
\Psi(G)\cdot \Psi(G')
\qquad ({^\forall}G,G'\in\mathbb{G}^{\bullet}).
\end{equation}
%
\end{definition}%

Let $(\mathbb{G}^{\bullet},R,\Psi)$ be a GA triple
and
$\mathbb{G}^{\circ}_d\subset \mathbb{G}_d^{\bullet}$ be the 
subset of connected graphs. 
By the finiteness  assumption \ref{finiteness}
of the G-set, the following sums over
$\mathbb{G}_d^{\bullet}$ and 
$\mathbb{G}^{\circ}_d$
are  
well-defined elements of $R$:
\begin{equation}\notag
\sum_{G\in \mathbb{G}^{\bullet}_d}
\frac{1}{\#{\aut}(G)}\Psi(G),
\qquad
\sum_{G\in \mathbb{G}^{\circ}_d}
\frac{1}{\#{\aut}(G)}\Psi(G).
\end{equation}
These are related by the relation:
\begin{prop}\label{expfml}
%
\begin{equation}\notag
1+\sum_{d\geq 1}
\sum_{G\in \mathbb{G}^{\bullet}_d}
\frac{1}{\#{\aut}(G)}\Psi(G)x^{d}
=\exp\Bigg[
\sum_{d\geq 1}\sum_{G\in \mathbb{G}^{\circ}_d}
\frac{1}{\#{\aut}(G)}\Psi(G)x^{d}
\Bigg]
\end{equation}
as formal power series in $x$.
\end{prop}
We call 
this the {\em exponential formula}
since it is
a variation of 
the exponential formula (\cite{stan}, chapter 5).

\begin{proof}%
Expanding the RHS, we obtain
\begin{equation}\notag
1+\sum_{d\geq 1}x^d
\Bigg(
\sum_{l\geq 1}
\sum_{\begin{subarray}{c}
\{(G_i,n_i)\}_{1\leq i\leq l};\\
G_i\in\mathbb{G}^{\circ},n_i\in\mathbb{N};\\
\sum_{i=1}^l n_i d(G_i)=d
\end{subarray}}
\underbrace{
\frac{1}{n_1!\cdots n_l!}
\prod_{i=1}^{l}\Big(\frac{\Psi(G_i)}{\#\aut(G_i)}\Big)^{n_i}
\Bigg)
}_{\ast}
.
\end{equation}
%
%
As we mentioned in remark \ref{gen},
each  $\{(G_1,n_1),\ldots,(G_l,n_l)\}$
uniquely corresponds to 
the graph 
\begin{equation}\notag
G=\bigcup_{i=1}^l 
\underbrace{G_i\cup\cdots\cup G_i}_{n_i\text{ times}},
\end{equation}
and $(\ast)$ is exactly the same as 
the contribution of $G$ to the LHS,
as we now explain.
The number of automorphisms of such $G$ is 
\begin{equation}\notag
\#{\aut}(G)=
\prod_{i=1}^l n_i!\big(\#{\aut}(G)\big)^{n_i}.
\end{equation}
The amplitude $\Psi(G)$ is 
the product of $\Psi(G_i)^{n_i}$ ($1\leq i\leq l$)
because of the multiplicativity of  $\Psi$.
Therefore
the contribution of $G$ to the LHS and 
$(\ast)$
are the same.
%
\end{proof}%

\section{Proof of Free energy}\label{pfe}
We give a proof of 
proposition \ref{free}
by using  the exponential formula.
All we have to do  is to define a
suitable GA triple.
We first construct a G-set.

For a combined forest $W$, we define $\overline{W}$ to be 
the graph obtained by forgetting all leaf indices.
Two  combined forests $W$ and $W'$ are {\em equivalent}
if $\overline{W}\cong\overline{W}'$.
The set of equivalence classes in
${\rm Comb}^{\bullet}(\Vec{\mu},\Vec{\nu},\Vec{\lambda})$ 
and
${\rm Comb}^{\circ}(\Vec{\mu},\Vec{\nu},\Vec{\lambda})$ 
are denoted by
$\overline{{\rm Comb}}^{\bullet}(\Vec{\mu},\Vec{\nu},\Vec{\lambda})$ 
and
$\overline{{\rm Comb}}^{\circ}(\Vec{\mu},\Vec{\nu},\Vec{\lambda})$.
Let us define
\begin{equation}\notag
\begin{split}
\overline{{\rm Comb}}^{\bullet}_{\gamma}(d)&=
\coprod_{|\Vec{d}|=d}
\coprod_{
\begin{subarray}{c}
(\Vec{\mu},\Vec{\nu},\Vec{\lambda})\\
\text{$r$-set of}\\
\text{degree $\Vec{d}$}
\end{subarray}}
\overline{{\rm Comb}}^{\bullet}_{\gamma}
(\Vec{\mu};\Vec{\nu};\Vec{\lambda}) \qquad (d\geq 1),
\\
\overline{{\rm Comb}}^{\bullet}_{\gamma}&=\coprod_{d\geq 1}
\overline{{\rm Comb}}_{\gamma}^{\bullet}(d).
\end{split}
\end{equation}
The set $\overline{{\rm Comb}}^{\bullet}_{\gamma}$
clearly satisfies the conditions 1,3,4,5 of the G-set.
It also satisfies the second condition, because
\begin{equation}\notag
\overline{{\rm Comb}}_{\gamma}^{\bullet}
(\vec{\mu},\vec{\nu},\vec{\lambda})
=\coprod_{k\geq 1}
\coprod_{\{(\Vec{\mu}^j,\Vec{\nu}^j,\Vec{\lambda}^j)\}_{1\leq j\leq k}
\in  D_k(\Vec{\mu},\Vec{\nu},\Vec{\lambda})}
\overline{{\rm Comb}}_{\gamma}^{\circ}
(\vec{\mu}^1,\vec{\nu}^1,\vec{\lambda}^1)\times
\ldots
\times
\overline{{\rm Comb}}_{\gamma}^{\circ}
(\vec{\mu}^k,\vec{\nu}^k,\vec{\lambda}^k)
\end{equation}
where $D_k$ is the set of all partitions of 
$(\Vec{\mu},\Vec{\nu},\Vec{\lambda})$
into $k$ parts 
$\{(\Vec{\mu}^j,\Vec{\nu}^j,\Vec{\lambda}^j)\}_{1\leq j\leq k}$
such that every part is an $r$-set.

Secondly, we define the amplitude $\Psi$ of
$G\in \overline{{\rm Comb}}_{\gamma}^{\bullet}
(\Vec{\mu},\Vec{\nu},\Vec{\lambda})$
so that it satisfies
\begin{equation}\notag
\frac{1}{\#\aut(G)}\Psi(G)
=
\sum_{\begin{subarray}{c}
W:\text{ combined forest}\\
\overline{W}=G
\end{subarray}
}
\frac{1}{z_{\Vec{\mu}}z_{\Vec{\nu}}z_{\Vec{\lambda}}}
{\mathcal H}(W) \,
\Vec{Q}^{\Vec{d}}
\end{equation}
where 
$\Vec{d}$ is the degree of the $r$-set
$(\vec{\mu},\vec{\nu},\vec{\lambda})$.
Solving this,
we have
\begin{equation}\notag
\Psi(G)=\#\aut(G)\cdot N(W)
\frac{1}
{z_{\Vec{\mu}}z_{\Vec{\nu}}z_{\Vec{\lambda}}}
{\mathcal H}(W)\Vec{Q}^{\Vec{d}}
\qquad 
(W \text{ such that } \overline{W}=G)
\end{equation}
where
$N(W)$ is 
the number of combined forests equivalent to
$W$.

It is easy to see that
$\Psi$ 
is a grade preserving map 
to
the ring
$R=\mathbb{Q}(t)[[\Vec{Q}]]$.
To check the multiplicativity of $\Psi$,
we need to know $N(W)$.
Since $N(W)$ is equal to the number of ways to assign
leaf indices to $\overline{W}$,
it consists of  two factors.
The one is the number of ways to distribute the 
leaf indices to VEV trees  in $\overline{W}$;
we define the partitions
$\lambda(T,T'),\mu(T),\nu(T)$ for
a VEV tree $T,T'$ in $W$ by:
\begin{equation}\notag
\begin{split}
\lambda(T,T')&=\{h(b)| \text{ $b$ is a bridge joining $T$, $T'$}\},
\\
\mu(T)&=\{c_v|
\text{ $v$ is a left leaf of $T$, not incident to a bridge}\},
\\
\nu(T)&=\{c_v|
\text{ $v$ is a right leaf of $T$, not incident to a bridge}\};
\end{split}
\end{equation}
then there are 
\begin{equation}\notag
\frac{\#\aut(\Vec{\mu})\,\#\aut(\Vec{\nu})\,\#\aut(\Vec{\lambda}) }
{|\aut(\overline{W})| 
\prod_{T\neq T'}\#\aut(\lambda(T,T'))
\prod_{T}
\#\aut(\mu(T))\,\#\aut(\nu(T))
}
\end{equation}
ways to distribute the leaf indices to VEV trees.
The other factor is the number of ways to assign the leaf indices
in each VEV tree.
It is the same as the number of (connected) 
combined forests
equivalent to $T$ which we write as $N'(T)$.

Therefore
\begin{equation}\notag
N(W)=
\frac{\#\aut(\Vec{\mu})\,\#\aut(\Vec{\nu})\,\#\aut(\Vec{\lambda}) }
{|\aut(\overline{W})| }
\underbrace{
\frac{1}{\prod_{T\neq T'}\#\aut(\lambda(T,T'))}
\prod_{T}\frac{N'(T)}
{\#\aut(\mu(T))\,\#\aut(\nu(T))
}
}_{\star}
\end{equation}
The factor $(\star)$ is multiplicative.
Non-multiplicative factors
cancel with those of 
$z_{\Vec{\mu}}$, $z_{\Vec{\nu}}$, $z_{\Vec{\lambda}}$
and
$\#\aut(\overline{W})$:
\begin{equation}\notag
\frac{\#\aut(G)\cdot N(W)}{z_{\Vec{\mu}}z_{\Vec{\nu}}z_{\Vec{\lambda}}}
=
{\prod_{i=1}^r \Big(
\prod_{j=1}^{l(\mu^i)} \mu^i_j \cdot
\prod_{j=1}^{l(\nu^i)} \nu^i_j \cdot 
\prod_{j=1}^{l(\lambda^i)} \lambda^i_j
}\Big)^{-1}
\times(\star).
\end{equation}
This is multiplicative and 
so are
${\mathcal H}(W)$ and $\Vec{Q}^{\Vec{d}}$.
Thus we finished the check of the multiplicativity of $\Psi$
and
$(\overline{{\rm Comb}}_{\gamma}^{\bullet},R,\Psi)$ 
is a GA triple.

With the GA triple, 
the partition function is expressed as
\begin{equation}\notag
{\mathcal Z}^{\gamma}(q;\Vec{Q})
=1+\sum_{G\in   \overline{{\rm Comb}}_{\gamma}^{\bullet}}
\frac{1}{\#\aut(G)}
\Psi(G).
\end{equation}
Therefore we apply
the exponential formula and obtain the free energy:
\begin{equation}\notag
\begin{split}
{\mathcal F}^{\gamma}(q;\Vec{Q})
&
=\sum_{G\in   \overline{{\rm Comb}}_{\gamma}^{\bullet}}
\frac{1}{\#\aut(G)}
\Psi(G)
\\
&=
\sum_{\Vec{d}\neq \Vec{0}}Q^{\Vec{d}}
\sum_{\begin{subarray}{c}
(\Vec{\mu},\Vec{\nu},\Vec{\lambda});\\
\text{$r$-set of}\\
\text{degree $\Vec{d}$}
\end{subarray}}
\frac{1}
{z_{\Vec{\lambda}}z_{\Vec{\mu}}z_{\Vec{\nu}}}
\sum_{W\in {\rm Comb}^{\circ}_{\gamma}
(\Vec{\mu},\Vec{\nu},\Vec{\lambda})}
{\mathcal H}(W).
\end{split}
\end{equation}
Since
${\mathcal F}_{\Vec{d}}^{\gamma}(q)$ is 
the coefficient of $\Vec{Q}^{\Vec{d}}$,
this implies the proposition \ref{free}.
\section{}\label{proofG4}
In the proof of proposition \ref{lastprop},
we have claimed the existence of a certain coupling
between VEV trees and bridges in 
a combined forest $W$.
If 
we make the simpler graph $\hat{W}$ by
contracting every VEV tree in 
$W$ to a vertex,
then
the coupling is the same as that of
vertices and edges in $\hat{W}$
(see, for example, figure \ref{excomb2} and (\ref{scd})).
Then the existence of such coupling 
follows 
if we apply
the next lemma
to $\hat{W}$.
%
\begin{lemma}\label{G4}%
Let $\Gamma$ be a connected 
graph without self-loops (a self loop is an edge that joins
a single vertex to itself). 
\begin{enumerate}
\item %
If the cycle rank $\beta(\Gamma)$ is zero, 
then for all $v\in V(\Gamma)$, there exists
a map $\varphi_v:E(\Gamma)\to V(\Gamma)$  satisfying
the following conditions;
$\varphi_v(e)$ is either of the incident vertices;
$|\varphi_v^{-1}(v')|=1$ for $v'\neq v$;
$|\varphi_v^{-1}(v)|=0$.  
We call $\varphi_v$ an {\em edge-map of $v$}. 
\item %
If the cycle rank $\beta(\Gamma)>1$,
then there exists
$\varphi: E(\Gamma)\to V(\Gamma)$ satisfying the following two conditions;
$\varphi(e)$ is either of the incident vertices;
$|\varphi^{-1}(v)|\geq 1$ for all $v\in V(\Gamma)$.
We  call $\varphi$ an {\em edge-map of $\Gamma$}. 
\end{enumerate}
\end{lemma}

\begin{proof}
Firstly,
we make a (finite) sequence $(\Gamma_i)_{i\geq 0}$
as follows.
\begin{enumerate}
\item 
$\Gamma_0=\Gamma$.
\item
Assume $\Gamma_i$ is given.
Then pick one vertex $v_i\in V(\Gamma_i)$.
(This choice is completely arbitrary.)
$\Gamma_{i+1}$ is the graph obtained by deleting
$v_i$ and all of its incident edges.
\item
Repeat the step 2 until we obtain the graph
consisting only of vertices.
Note that this process ends after at most
$\#V(\Gamma)$ steps.
\end{enumerate}
We write the sequence as follows.
\begin{equation}\notag
\Gamma=\Gamma_0\stackrel{}{\longrightarrow}
\Gamma_1\stackrel{}{\longrightarrow}
\cdots
\stackrel{}{\longrightarrow}
\Gamma_n=(\text{vertices})
\end{equation}
As examples, consider the figures in (\ref{fst})(\ref{scd}).
We can take the following sequences:
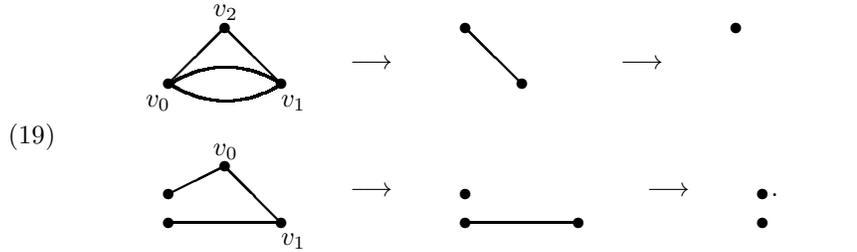
\newcommand{\fstg}{
\unitlength .15cm
\begin{picture}(15,10)(-5,-4)
\thicklines
\put(0,0){\line(1,1){5}}
\put(10,0){\line(-1,1){5}}
\qbezier(0,0)(5,3)(10,0)
\qbezier(0,0)(5,-3)(10,0)
\put(0,0){\circle*{1}}
\put(10,0){\circle*{1}}
\put(5,5){\circle*{1}}
\put(-2,-2){$v_0$}
\put(10,-2){$v_1$}
\put(4,6){$v_2$}
\end{picture}
}
\newcommand{\fstgr}{
\unitlength .15cm
\begin{picture}(15,10)(0,-4)
\thicklines
%
\put(10,0){\line(-1,1){5}}
\put(10,0){\circle*{1}}
\put(5,5){\circle*{1}}
\end{picture}
}
\newcommand{\fstgra}{
\unitlength .15cm
\begin{picture}(15,10)(0,-4)
\thicklines
%
\put(5,5){\circle*{1}}
\end{picture}
}
\newcommand{\scdg}{
\unitlength .15cm
\begin{picture}(15,10)(-5,-3)
\thicklines
\put(0,2.5){\line(2,1){5}}
\put(10,0){\line(-1,1){5}}
\put(0,0){\line(1,0){10}}
\put(0,0){\circle*{1}}
\put(0,2.5){\circle*{1}}
\put(10,0){\circle*{1}}
\put(5,5){\circle*{1}}
\put(10,-2){$v_1$}
\put(4,6){$v_0$}
\end{picture}
}
\newcommand{\scdgr}{
\unitlength .15cm
\begin{picture}(15,10)(-5,-3)
\thicklines
\put(0,0){\line(1,0){10}}
\put(0,0){\circle*{1}}
\put(0,2.5){\circle*{1}}
\put(10,0){\circle*{1}}
\end{picture}
}
\newcommand{\scdgra}{
\unitlength .15cm
\begin{picture}(5,10)(-5,-3)
\thicklines
\put(0,0){\circle*{1}}
\put(0,2.5){\circle*{1}}
\end{picture}
}
\begin{equation}\label{seq}
\begin{split}
&
\raisebox{-.8cm}{\fstg}\qquad\longrightarrow
\raisebox{-.8cm}{\fstgr}\quad\longrightarrow
\raisebox{-.8cm}{\fstgra},
\\
&
\raisebox{-.8cm}{\scdg}\qquad\longrightarrow
\raisebox{-.8cm}{\scdgr}\quad\quad\longrightarrow
\raisebox{-.8cm}{\scdgra}.
\end{split}
\end{equation}

Now, by induction, we will construct a map 
\begin{equation}\notag
\varphi_i:E(\Gamma_i)\to V(\Gamma_i)
\end{equation}
such that $\varphi_i(e)$ is one of the endpoints of an edge $e$.
\begin{enumerate}
\item 
We define $\varphi_n$ to be the empty map
since there is no edges in $\Gamma_n$.
\item
Next assume that 
$\varphi_{i+1}$ is given.
If an edge $e\in E(\Gamma_i)$ is not incident to $v_i$,
it has the corresponding edge $e'$ in
in $\Gamma_{i+1}$. We define
\begin{equation}\notag
\varphi_i(e)=\varphi_{i+1}(e')
\end{equation}
If an edge $e\in E(\Gamma_i)$ is incident to $v_i$,
%
we define
\begin{equation}\notag
\begin{cases}
\varphi_i(e)=v_i
&\text{ if }\varphi_{i+1}^{-1}(v)\neq \emptyset\\
\varphi_i(e)=v 
&\text{ if }\varphi_{i+1}^{-1}(v)=\emptyset
\end{cases}
\end{equation}%
where $v$ is the other endpoint of $e$.
\end{enumerate}
Thus we obtain
$\varphi_0:E(\Gamma)\to V(\Gamma)$.
If $\Gamma$ has the cycle rank $>0$,
$\varphi_0$ is its edge-map.
If the cycle rank is zero,
$\varphi_0$ is an edge-map of $v_0$.
We make an edge-map of other vertex $v\in V(\Gamma)$ as follows.
We take the path
$v=v_{(0)}\to v_{(1)}\to\cdots$ 
where
$e_{(i)}=\varphi_{v_i}^{-1}(v_{(i)})$
and $v_{(i+1)}$ is the other endpoint of $e_{(i)}$:
\newcommand{\mypath}{  
\unitlength .2cm
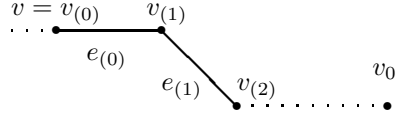
\begin{picture}(15,10)(-5,-1)
\thicklines
\multiput(-5,5)(1,0){5}{\line(1,0){.1}}
\put(-2,5){\line(1,0){7}}
\put(5,5){\line(1,-1){5}}
\multiput(10,0)(1,0){10}{\line(1,0){.1}}
\put(-2,5){\circle*{.5}}
\put(5,5){\circle*{.5}}
\put(10,0){\circle*{.5}}
\put(20,0){\circle*{.5}}
\put(-5,6){$v=v_{(0)}$}
\put(4,6){$v_{(1)}$}
\put(10,1){$v_{(2)}$}
\put(0,3){$e_{(0)}$}
\put(5,1){$e_{(1)}$}
\put(19,2){$v_0$}
\end{picture}
}
\begin{equation}\notag
\mypath
\end{equation}
Since the cycle-rank is zero, this path arrives at $v_0$
at some time.
Define
\begin{equation}\notag
\varphi_v(e)=\begin{cases}
\varphi_0(e)&\text{ if $e$ is not on the path}\\
v_{(i+1)}& \text{ if $e=e_{(i)}$}.
\end{cases}
\end{equation}
Then $\varphi_v$ is an edge-map of $v$.
This completes the proof.
\end{proof}
The examples of coupling given 
after (\ref{fst}) and (\ref{scd})
were constructed in this way from (\ref{seq}).

%

\end{document}

%% file: graphexp2.tex

\newcommand{\oneone}[1]{{  
\unitlength .1cm
\begin{picture}(15,12)(-5,-3)
\thicklines
\put(0,5){\circle*{1}}
\put(10,5){\circle*{1}}
\put(5,0){\circle*{1}}
\put(5,0){\line(-1,1){5}}
\put(5,0){\line(1,1){5}}
\put(-1.5,6){$#1$}
\put(8,6){$-#1$}
\end{picture}
}}

\newcommand{\onetwothree}[3]{{  
\unitlength .1cm
\begin{picture}(15,17)(-15,-3)
\thicklines
\put(0,0){\circle*{1}}
\put(5,5){\circle*{1}}
\put(10,10){\circle*{1}}
\put(0,10){\circle*{1}}
\put(-10,10){\circle*{1}}
\put(0,0){\line(1,1){10}}
\put(0,0){\line(-1,1){10}}
\put(5,5){\line(-1,1){5}}
\put(8,11){$-#3$}
\put(-11,11){$#1$}
\put(-1,11){$#2$}
\end{picture}
}}
\newcommand{\onetofourth}[2]{{  
\unitlength .1cm
\begin{picture}(20,22)(-20,-3)
\thicklines
\put(0,0){\circle*{1}}
\put(5,5){\circle*{1}}
\put(15,15){\circle*{1}}
\put(-15,15){\circle*{1}}
\put(0,10){\circle*{1}}
\put(-5,15){\circle*{1}}
\put(5,15){\circle*{1}}
\put(0,0){\line(1,1){15}}
\put(0,0){\line(-1,1){15}}
\put(5,5){\line(-1,1){10}}
\put(0,10){\line(1,1){5}}
\put(-17,16){$#1$}
\put(-7,16){$#2$}
\put(3,16){$-#2$}
\put(13,16){$-#1$}
\end{picture}
}}
\newcommand{\twooneoneichi}{{  
\unitlength .1cm
\begin{picture}(20,17)(-15,-3)
\thicklines
\put(0,0){\circle*{1}}
\put(-10,10){\circle*{1}}
\put(0,10){\circle*{1}}
\put(10,10){\circle*{1}}
\put(20,10){\circle*{1}}
\put(10,0){\circle*{1}}
\put(0,0){\line(-1,1){10}}
\put(0,0){\line(1,1){10}}
\put(10,0){\line(1,1){10}}
\put(10,0){\line(-1,1){4}}
\put(0,10){\line(1,-1){4}}
\put(-12,11){$1$}
\put(-2,11){$1$}
\put(8,11){$-1$}
\put(18,11){$-1$}
\end{picture}}}
\newcommand{\twooneoneni}[2]{{  
\unitlength .1cm
\begin{picture}(20,22)(-20,-3)
\thicklines
\put(0,0){\circle*{1}}
\put(-15,15){\circle*{1}}
\put(-5,15){\circle*{1}}
\put(5,15){\circle*{1}}
\put(15,15){\circle*{1}}
\put(0,10){\circle*{1}}
\put(0,0){\line(1,1){15}}
\put(0,0){\line(-1,1){15}}
\put(0,10){\line(-1,1){5}}
\put(0,10){\line(1,1){5}}
\put(-17,16){$#1$}
\put(-7,16){$#2$}
\put(3,16){$-#2$}
\put(13,16){$-#1$}
\end{picture}}}
\newcommand{\oneoneonethree}{{  
\unitlength .1cm
\begin{picture}(20,22)(-20,-3)
\thicklines
\put(0,0){\circle*{1}}
\put(-15,15){\circle*{1}}
\put(-5,15){\circle*{1}}
\put(5,15){\circle*{1}}
\put(15,15){\circle*{1}}
\put(5,5){\circle*{1}}
\put(10,10){\circle*{1}}
\put(0,0){\line(1,1){15}}
\put(0,0){\line(-1,1){15}}
\put(5,5){\line(-1,1){10}}
\put(10,10){\line(-1,1){5}}
\put(-17,16){$1$}
\put(-7,16){$1$}
\put(3,16){$1$}
\put(12,16){$-3$}
\end{picture}}}

\newcommand{\oneoneonetwoone}{{  
\unitlength .1cm
\begin{picture}(30,20)(-10,-3)
\thicklines
\put(0,0){\circle*{1}}
\put(-15,15){\circle*{1}}
\put(-5,15){\circle*{1}}
\put(5,15){\circle*{1}}
\put(15,15){\circle*{1}}
\put(25,15){\circle*{1}}
\put(5,5){\circle*{1}}
\put(15,5){\circle*{1}}
\put(0,0){\line(1,1){15}}
\put(0,0){\line(-1,1){15}}
\put(5,5){\line(-1,1){10}}
\put(15,5){\line(-1,1){4}}
\put(15,5){\line(1,1){10}}
\put(5,15){\line(1,-1){4}}
\put(-17,16){$1$}
\put(-7,16){$1$}
\put(3,16){$1$}
\put(13,16){$-2$}
\put(23,16){$-1$}
\end{picture}}}
\newcommand{\oneoneonetwooneni}{{  
\unitlength .1cm
\begin{picture}(30,20)(-10,-3)
\thicklines
\put(0,0){\circle*{1}}
\put(-15,15){\circle*{1}}
\put(-5,15){\circle*{1}}
\put(5,15){\circle*{1}}
\put(15,15){\circle*{1}}
\put(25,15){\circle*{1}}
\put(10,10){\circle*{1}}
\put(10,0){\circle*{1}}
\put(0,0){\line(1,1){15}}
\put(0,0){\line(-1,1){15}}
\put(10,10){\line(-1,1){5}}
\put(10,0){\line(1,1){15}}
\put(10,0){\line(-1,1){4}}
\put(-5,15){\line(1,-1){9}}
\put(-17,16){$1$}
\put(-7,16){$1$}
\put(3,16){$1$}
\put(13,16){$-2$}
\put(23,16){$-1$}

\end{picture}}}
\newcommand{\oneoneonetwoonesan}{{  
\unitlength .1cm
\begin{picture}(25,27)(-15,-3)
\thicklines
\put(0,0){\circle*{1}}
\put(-20,20){\circle*{1}}
\put(20,20){\circle*{1}}
\put(-10,20){\circle*{1}}
\put(0,20){\circle*{1}}
\put(10,20){\circle*{1}}
\put(0,10){\circle*{1}}
\put(5,15){\circle*{1}}
\put(0,0){\line(1,1){20}}
\put(0,0){\line(-1,1){20}}
\put(0,10){\line(1,1){10}}
\put(0,10){\line(-1,1){10}}
\put(5,15){\line(-1,1){5}}
\put(-22,21){$1$}
\put(-12,21){$1$}
\put(-2,21){$1$}
\put(7,21){$-2$}
\put(17,21){$-1$}
\end{picture}}}
\newcommand{\oneoneonetwooneyon}{{  
\unitlength .1cm
\begin{picture}(25,27)(-15,-3)
\thicklines
\put(0,0){\circle*{1}}
\put(-20,20){\circle*{1}}
\put(20,20){\circle*{1}}
\put(-10,20){\circle*{1}}
\put(0,20){\circle*{1}}
\put(10,20){\circle*{1}}
\put(0,10){\circle*{1}}
\put(5,15){\circle*{1}}
\put(5,5){\circle*{1}}
\put(0,0){\line(1,1){20}}
\put(0,0){\line(-1,1){20}}
\put(0,10){\line(1,1){10}}
\put(0,10){\line(-1,1){10}}
\put(5,15){\line(-1,1){5}}
\put(0,10){\line(1,-1){5}}
\put(-22,21){$1$}
\put(-12,21){$1$}
\put(-2,21){$1$}
\put(7,21){$-2$}
\put(17,21){$-1$}
\end{picture}}}
\newcommand{\oneonethirdichi}{{  
\unitlength .1cm
\begin{picture}(30,32)(-20,-3)
\thicklines
\put(0,0){\circle*{1}}
\put(-25,25){\circle*{1}}
\put(-15,25){\circle*{1}}
\put(-5,25){\circle*{1}}
\put(5,25){\circle*{1}}
\put(15,25){\circle*{1}}
\put(25,25){\circle*{1}}
\put(5,5){\circle*{1}}
\put(10,10){\circle*{1}}
\put(-5,15){\circle*{1}}
\put(5,15){\circle*{1}}
\put(0,0){\line(1,1){25}}
\put(0,0){\line(-1,1){25}}
\put(5,5){\line(-1,1){20}}
\put(10,10){\line(-1,1){9}}
\put(-5,15){\line(1,1){10}}
\put(5,15){\line(1,1){10}}
\put(-5,25){\line(1,-1){4}}
\put(-27,26){$1$}
\put(-17,26){$1$}
\put(-7,26){$1$}
\put(2,26){$-1$}
\put(12,26){$-1$}
\put(22,26){$-1$}
\end{picture}}}

\newcommand{\oneonethirdni}{{  
\unitlength .1cm
\begin{picture}(30,32)(-20,-3)
\thicklines
\put(0,0){\circle*{1}}
\put(-25,25){\circle*{1}}
\put(-15,25){\circle*{1}}
\put(-5,25){\circle*{1}}
\put(5,25){\circle*{1}}
\put(15,25){\circle*{1}}
\put(25,25){\circle*{1}}
\put(5,5){\circle*{1}}
\put(0,20){\circle*{1}}
\put(10,10){\circle*{1}}
\put(0,10){\circle*{1}}
\put(0,0){\line(1,1){25}}
\put(0,0){\line(-1,1){25}}
\put(5,5){\line(-1,1){20}}
\put(-5,25){\line(1,-1){9}}
\put(10,10){\line(-1,1){4}}
\put(0,10){\line(1,1){15}}
\put(0,20){\line(1,1){5}}

\put(-27,26){$1$}
\put(-17,26){$1$}
\put(-7,26){$1$}
\put(2,26){$-1$}
\put(12,26){$-1$}
\put(22,26){$-1$}

\end{picture}}}

\newcommand{\oneonethirdsan}{{  
\unitlength .1cm
\begin{picture}(30,32)(-20,-3)
\thicklines
\put(0,0){\circle*{1}}
\put(-25,25){\circle*{1}}
\put(-15,25){\circle*{1}}
\put(-5,25){\circle*{1}}
\put(5,25){\circle*{1}}
\put(15,25){\circle*{1}}
\put(25,25){\circle*{1}}
\put(5,5){\circle*{1}}
\put(5,15){\circle*{1}}
\put(0,10){\circle*{1}}
\put(0,20){\circle*{1}}
\put(0,0){\line(1,1){25}}
\put(0,0){\line(-1,1){25}}
\put(5,5){\line(-1,1){20}}
\put(0,10){\line(1,1){15}}
\put(5,15){\line(-1,1){10}}
\put(0,20){\line(1,1){5}}
\put(-27,26){$1$}
\put(-17,26){$1$}
\put(-7,26){$1$}
\put(2,26){$-1$}
\put(12,26){$-1$}
\put(22,26){$-1$}

\end{picture}}}
\newcommand{\treea}{{  
\unitlength .1cm
\begin{picture}(17,17)(-17,-2)
\thicklines
\put(0,0){\line(1,1){15}}
\put(0,0){\line(-1,1){15}}
\put(5,5){\line(-1,1){10}}
\put(0,10){\line(1,1){5}}
\put(0,0){\circle*{1}}
\put(5,5){\circle*{1}}
\put(15,15){\circle*{1}}
\put(-15,15){\circle*{1}}
\put(0,10){\circle*{1}}
\put(-5,15){\circle*{1}}
\put(5,15){\circle*{1}}
\end{picture}
}}
\newcommand{\treeb}{{  
\unitlength .1cm
\begin{picture}(17,12)(-12,-7)
\thicklines
\put(0,0){\line(-1,1){10}}
\put(0,0){\line(1,1){10}}
\put(10,0){\line(1,1){10}}
\put(10,0){\line(-1,1){4}}
\put(0,10){\line(1,-1){4}}
\put(0,0){\circle*{1}}
\put(-10,10){\circle*{1}}
\put(0,10){\circle*{1}}
\put(10,10){\circle*{1}}
\put(20,10){\circle*{1}}
\put(10,0){\circle*{1}}
\end{picture}}}
\newcommand{\treec}{{  
\unitlength .1cm
\begin{picture}(17,17)(-17,-2)
\thicklines
\put(0,0){\line(1,1){15}}
\put(0,0){\line(-1,1){15}}
\put(0,10){\line(-1,1){5}}
\put(0,10){\line(1,1){5}}
\put(0,0){\circle*{1}}
\put(-15,15){\circle*{1}}
\put(-5,15){\circle*{1}}
\put(5,15){\circle*{1}}
\put(15,15){\circle*{1}}
\put(0,10){\circle*{1}}
\end{picture}}}
\newcommand{\treed}{{  
\unitlength .1cm
\begin{picture}(12,7)(-2,-2)
\thicklines
\put(5,0){\line(-1,1){5}}
\put(5,0){\line(1,1){5}}
\put(0,5){\circle*{1}}
\put(10,5){\circle*{1}}
\put(5,0){\circle*{1}}
\end{picture}
}}

\newcommand{\comb}[2]{{
\unitlength .1cm
\begin{picture}(50,45)(-50,-15)

\put(-43,-12){#1}
\put(7,-12){#2}
\put(-8,18){\treed}
\put(24,-16){$2$}
\put(-26,-16){$1$}
\put(-1,16){$3$}

\qbezier(-10,5)(0,15)(10,5)
\qbezier(-20,5)(0,25)(20,5)
\qbezier(-5,25)(15,40)(40,5)
\qbezier(5,25)(-15,40)(-40,5)

\end{picture}
}}

\newcommand{\ancomb}[2]{{
\unitlength .1cm
\begin{picture}(50,45)(-50,-15)
\put(-43,-12){#1}
\put(7,-12){#2}
\put(-8,18){\treed}
\put(24,-16){$2$}
\put(-26,-16){$1$}
\put(-1,16){$3$}
\qbezier(-10,5)(0,15)(10,5)
\qbezier(-5,25)(15,40)(40,5)
\qbezier(5,25)(-15,40)(-40,5)
\put(-30,25){$b_3$}
\put(27,25){$b_2$}
\put(0,5){$b_1$}
\put(-33,-10){$T_1$}
\put(-20,-10){$T_2$}
\put(17,-10){$T_3$}
\put(7,20){$T_4$}
\end{picture}
}}

\newcommand{\poichi}{{  
\unitlength .1cm
\begin{picture}(30,30)(-30,-3)
\thicklines
\put(0,0){\circle*{1}}
\put(-25,25){\circle*{1}}
\put(-15,25){\circle*{1}}
\put(-5,25){\circle*{1}}
\put(5,25){\circle*{1}}
\put(15,25){\circle*{1}}
\put(25,25){\circle*{1}}
\put(10,20){\circle*{1}}
\put(-10,20){\circle*{1}}
\put(5,5){\circle*{1}}
\put(15,15){\circle*{1}}
\put(0,0){\line(-1,1){25}}
\put(0,0){\line(1,1){25}}
\put(5,5){\line(-1,1){20}}
\put(15,15){\line(-1,1){10}}
\put(10,20){\line(1,1){5}}
\put(-10,20){\line(1,1){5}}
\put(-26,27){$a$}
\put(-16,27){$b$}
\put(-6,27){$c$}
\put(4,27){$d$}
\put(14,27){$e$}
\put(24,27){$f$}
\end{picture}
}}

\newcommand{\poni}{{  
\unitlength .1cm
\begin{picture}(30,35)(-30,-3)
\thicklines
\put(0,0){\circle*{1}}
\put(-25,25){\circle*{1}}
\put(-15,25){\circle*{1}}
\put(-5,25){\circle*{1}}
\put(5,25){\circle*{1}}
\put(15,25){\circle*{1}}
\put(25,25){\circle*{1}}
\put(10,20){\circle*{1}}
\put(-10,20){\circle*{1}}
\put(5,5){\circle*{1}}
\put(15,15){\circle*{1}}
\put(0,0){\line(-1,1){25}}
\put(0,0){\line(1,1){25}}
\put(5,5){\line(-1,1){20}}
\put(15,15){\line(-1,1){10}}
\put(10,20){\line(1,1){5}}
\put(-10,20){\line(1,1){5}}
\put(-26,27){$a$}
\put(-26,30){$1$}
\put(-16,27){$b$}
\put(-16,30){$2$}
\put(-6,27){$c$}
\put(-6,30){$3$}
\put(4,27){$d$}
\put(4,30){$4$}
\put(14,27){$e$}
\put(14,30){$5$}
\put(24,27){$f$}
\put(24,30){$6$}
\end{picture}
}}

\newcommand{\poyon}{{  
\unitlength .1cm
\begin{picture}(30,35)(-30,-3)
\thicklines
\put(0,0){\circle*{1}}
\put(-25,25){\circle*{1}}
\put(-15,25){\circle*{1}}
\put(-5,25){\circle*{1}}
\put(5,25){\circle*{1}}
\put(15,25){\circle*{1}}
\put(25,25){\circle*{1}}
\put(0,20){\circle*{1}}
\put(0,10){\circle*{1}}
\put(5,5){\circle*{1}}
\put(10,10){\circle*{1}}
\put(0,0){\line(-1,1){25}}
\put(0,0){\line(1,1){25}}
\put(5,5){\line(-1,1){20}}
\put(10,10){\line(-1,1){15}}
\put(0,20){\line(1,1){5}}
\put(0,10){\line(1,1){15}}
\put(-26,27){$a$}
\put(-26,30){$1$}
\put(-16,27){$b$}
\put(-16,30){$2$}
\put(-6,27){$d$}
\put(-6,30){$3$}
\put(4,27){$e$}
\put(4,30){$4$}
\put(14,27){$c$}
\put(14,30){$5$}
\put(24,27){$f$}
\put(24,30){$6$}
\end{picture}
}}

\newcommand{\posan}{{  
\unitlength .1cm
\begin{picture}(30,35)(-30,-3)
\thicklines
\put(0,0){\circle*{1}}
\put(-25,25){\circle*{1}}
\put(-15,25){\circle*{1}}
\put(-5,25){\circle*{1}}
\put(5,25){\circle*{1}}
\put(15,25){\circle*{1}}
\put(25,25){\circle*{1}}
\put(5,15){\circle*{1}}
\put(-5,15){\circle*{1}}
\put(5,5){\circle*{1}}
\put(10,10){\circle*{1}}
\put(0,0){\line(-1,1){25}}
\put(0,0){\line(1,1){25}}
\put(5,5){\line(-1,1){20}}
\put(10,10){\line(-1,1){15}}
\put(-5,15){\line(1,1){10}}
\put(5,15){\line(1,1){10}}
\put(-26,27){$a$}
\put(-26,30){$1$}
\put(-16,27){$b$}
\put(-16,30){$2$}
\put(-6,27){$d$}
\put(-6,30){$3$}
\put(4,27){$c$}
\put(4,30){$4$}
\put(14,27){$e$}
\put(14,30){$5$}
\put(24,27){$f$}
\put(24,30){$6$}
\end{picture}
}}
